\documentclass[letterpaper,12pt]{article}

\usepackage{color}
\usepackage{amssymb}
\usepackage{amsmath}
\usepackage{amscd}
\usepackage{mathtools}
\usepackage[all]{xy}

\usepackage[backref=page]{hyperref}
\hypersetup{colorlinks=true,linkcolor=blue,anchorcolor=blue,citecolor=blue}

\newtheorem{theo}{Theorem}[section]
\newtheorem{prop}[theo]{Proposition}
\newtheorem{lem}[theo]{Lemma}
\newtheorem{cor}[theo]{Corollary}
\newtheorem{defi}[theo]{Definition}
\newtheorem{rem}[theo]{Remark}

\def \K{{\mathcal K}}

\def \Br {{\rm{Br}}}

\def \si {{\sigma}}
\def \Ga {{\Gamma}}
\def \R {{\mathbb{R}}}
\def \Pic {{\rm {Pic}}}
\def \Div {{\rm {Div}}}
\def \div {{\rm{div}}}
\def \Gal {{\rm{Gal}}}
\def \Ker {{\rm{Ker}}}
\def \Im {{\rm {Im}}}
\def \Frob{{\rm{Frob}}}

\def \A{{\mathbb A}}
\def \P{{\mathbb P}}
\def \Spec {{\rm{Spec}}}
\def \dim {{\rm{dim}}}
\def \Hom {{\rm {Hom}}}

\def \Aut{{\rm Aut}}

\def\ov{\overline}

\def \Z {{\mathbb Z}}
\def \Q {{\mathbb Q}}
\def \F {{\mathbb F}}

\def \ZZ {{\rm Z}}
\def \Tr {{\rm{Tr}}}

\def\G{{\mathbb G}}

\def\C{{\mathbb C}}
\def\B{{\cal B}}

\def\lra{\longrightarrow}

\def\H{{\rm H}}

\def\Tr{{\rm Tr}}

\def\K{{\cal K}}

\def\O{{\cal O}}

\def\si{\sigma}

\def\Ga{\Gamma}

\def\e{\varepsilon}
\def\et{{\rm{\acute et}}}

\newcommand{\bthe}{\begin{theo}}
\newcommand{\ble}{\begin{lem}}
\newcommand{\bpr}{\begin{prop}}
\newcommand{\bco}{\begin{cor}}
\newcommand{\bde}{\begin{defi}}
\newcommand{\brem}{\begin{rem}}
\newcommand{\ethe}{\end{theo}}
\newcommand{\ele}{\end{lem}}
\newcommand{\epr}{\end{prop}}
\newcommand{\eco}{\end{cor}}
\newcommand{\ede}{\end{defi}}
\newcommand{\erem}{\end{rem}}

\title{Cohomology and the Brauer groups of diagonal surfaces}
\author{Dami\'an Gvirtz and Alexei N. Skorobogatov}
\date{\today}
\begin{document}
\baselineskip=15pt
\maketitle

\begin{flushright} {\em in memory of Sir Peter Swinnerton-Dyer}
\end{flushright}

\begin{abstract} \noindent  
We present a method for calculating the Brauer group of a surface 
given by a diagonal equation in projective space.
For diagonal quartic surfaces with coefficients in $\Q$
we determine the Brauer groups over $\Q$ and $\Q(i)$. 
\end{abstract}

\tableofcontents

\section*{Introduction}

Diagonal equations occupy a special place in number theory:
due to separation of variables many questions about such equations can be
reduced to questions about equations in one variable. This provided Andr\'e Weil with a key
example for the Weil conjectures \cite{W1, W2}. Diagonal hypersurfaces
are a particular case of varieties dominated by products of curves (DPC).
The Tate conjecture is known for DPC varieties defined over 
fields finitely generated over $\Q$, see \cite[(5.5)]{Tate}.
Using the results of \cite{SZ14} (based on the work of Faltings) it is easy to show that for DPC varieties 
defined over such fields the Galois-invariant subgroup
of the geometric Brauer group is finite. 

In this paper we use diagonal surfaces to test new techniques for 
computing the Brauer group, a task motivated by the Brauer--Manin obstruction 
to the local-to-global principle for rational points over number fields (see \cite{CTS21}
for a general introduction).
To put this into perspective, in characteristic zero, the Brauer group of any quadric
$X$ over any field $k$ is trivial, that is, equal to the image $\Br_0(X)$ 
of the canonical map $\Br(k)\to\Br(X)$.
The Brauer group of a cubic surface $X$ can be non-trivial, but it is equal to
the {\em algebraic Brauer group}, that is, to 
the kernel $\Br_1(X)$ of the canonical map $\Br(X)\to \Br(\ov X)$,
where $\ov X=X\times_k\bar k$ for an algebraic closure $\bar k$ of $k$.
This simplifies matters considerably, because $\Br_1(X)/\Br_0(X)$ is a subgroup of 
$\H^1(k,\Pic(\ov X))$ and is equal to this group if $k$ is a number field.
For diagonal cubic surfaces over $\Q$ all elements
of $\Br(X)$ modulo $\Br_0(X)$ and the ensuing Brauer--Manin conditions were determined in \cite{CTKS}.

The arithmetic of diagonal quartic surfaces is a recurrent theme in the work of Peter Swinnerton-Dyer,
from his first paper \cite{SD43}, revisited in \cite{SD68, PSD, SD00, SD14}, 
until his last published paper \cite{SD16}.
The algebraic Brauer groups $\Br_1(X)$ of diagonal quartic surfaces over arbitrary fields of
characteristic zero were classified by Bright in his thesis \cite{Bright, Bright2}, though an explicit form of 
the elements of this group is not known in general, cf.~\cite{Bright3}.
Ieronymou \cite[Thm.~3.1]{I} constructed generators of $\Br(\ov X)[2]$
and showed \cite[Prop.~4.9]{I} that they come from $4$-torsion elements of
the Brauer group of the quartic Fermat surface over $\Q(i,\sqrt[4]{2})$
(in fact, they come from $2$-torsion elements, see Corollary \ref{c11} below).
He then obtained sufficient conditions on the diagonal quartic surface $X$
over $\Q$ under which the 2-torsion subgroup
of $\Br(X)/\Br_1(X)$ is trivial \cite[Thm.~5.2]{I}. This, together with Mizukami's isomorphism between
the quartic Fermat surface over $\Q(i,\sqrt{2})$ and a certain Kummer surface, 
was used in \cite{ISZ} to deduce an explicit upper bound for the size of $\Br(X)/\Br_0(X)$. 
Finally, the odd order torsion subgroup of $\Br(X)$ was determined in \cite{IS} 
using a different method; 
the elements of this group have an explicit description, enabling one to compute the associated
Brauer--Manin conditions. The complete determination of the 2-primary torsion subgroup of $\Br(X)$
could not be done by previously known methods. We now describe a new approach
that allowed us to achieve this.

Let $k$ be a number field and let $\Ga=\Gal(\bar k/k)$.
For a smooth and projective {\em surface} $X$ over $k$
such that $\Pic(\ov X)$ is torsion-free and $\Br(\ov X)^\Ga$ is finite,
the determination of $\Br(X)/\Br_0(X)$ involves the following steps. 
(These conditions hold, for example, when $X$ is a K3 surface \cite{SZ08} or
a diagonal surface in $\P^3_k$, see \cite[XI.1.8]{SGA} and Proposition \ref{19nov} below.)

\smallskip

(a) Determine the action of $\Ga$
on the {\em geometric Brauer group} $\Br(\ov X)$, hence compute $\Br(\ov X)^\Ga$.
For a K3 surface 
with complex multiplication this involves identifying the Gr\"ossencharakter 
which describes the action of $\Ga$ on the Tate module of $\Br(\ov X)$. 
For diagonal quartic surfaces, this has essentially been done by 
Pinch and Swinnerton-Dyer in \cite{PSD} using work of Weil \cite{W2}.

(b) Determine the image of $\Br(X)\to \Br(\ov X)^\Ga$, called the {\em transcendental Brauer group}. 
The crucial fact is that when $\Pic(\ov X)$ is torsion-free, 
the 2-term complex of $\Ga$-modules 
\begin{equation}\Hom(\Pic(X_\C),\Z)\lra T(X_\C)\otimes\Q/\Z\label{co1}
\end{equation}
is quasi-isomorphic to the truncated complex $\tau_{[1,2]}{\bf R}p_*\G_{m,X}$,
where $p\colon X\to \Spec(k)$ is the structure morphism. Here $T(X_\C)$ is the {\em
transcendental lattice}, and the map in (\ref{co1}) comes from the cup-product
bilinear form on $\H^2(X_\C,\Z(1))$.

(c) The group $\Br_1(X)/\Br_0(X)\cong\H^1(k,\Pic(\ov X))$ is finite; 
it can be determined if we know an integral basis
of $\Pic(\ov X)$ and the action of $\Ga$ on it. The group $\Br(X)/\Br_0(X)$, an extension of 
finite abelian groups $\Br(X)/\Br_1(X)$ by $\Br_1(X)/\Br_0(X)$,
is computed as the first Galois hypercohomology group of the complex (\ref{co1}).

\smallskip

The methods in (b) and (c) are used for the first time in this paper. We apply them to
advance this program for diagonal surfaces of arbitrary degree $d$ in $\P^3_k$, using
the description of integral homology of the Fermat hypersurface $F$ of degree $d$ over $\C$
obtained by Pham \cite{Ph}, see also \cite{L, ABB, DS}.

Pham considers a natural ``vanishing cycle" in the complement $U_\C$ to a 
coordinate hyperplane section of $F_\C$
and shows that its orbit under the action of the group of diagonal automorphisms $G=(\mu_d)^3$ generates
the primitive subgroup of $\H_2(F_\C,\Z)$ (the kernel of the intersection product with the hyperplane
section). He also computes the intersection product in the resulting integral basis.
The group $\H^2(F_\C,\C)$ decomposes into a direct sum of 
1-dimensional eigenspaces of $G$, compatibly with the Hodge decomposition (which was determined by
Shioda \cite{Sh79} using work of Griffiths).
The Galois representation in the \'etale cohomology of $F$ goes back to
Weil \cite{W1, W2}, and was studied by Katz, Shioda, and Ulmer, among others, see \cite{U}. 
Combining these data, the computation of the Brauer group of diagonal surfaces
of any degree $d$ over the cyclotomic field $\Q(\mu_d)$ can in principle be done by a computer.
We work out precise details in the case of diagonal quartic surfaces with coefficients in $\Q$
over the ground field $\Q(i)$,
which are then used to deal with the case of the ground field~$\Q$. 

To state our results, we call two diagonal quartic forms
$\sum_{i=0}^3 a_ix_i^4$, where $a_i\in K^\times$ for $i=0,1,2,3$,
{\em equivalent} if one is obtained from another by permuting the variables $x_0,x_1,x_2,x_3$,
multiplying the coefficients $a_i$ by fourth powers in $K^\times$, and multiplying all four coefficients
by a common multiple in $K^\times$. A diagonal quartic surface is the surface given by
a diagonal quartic form in $\P^3_K$; surfaces given by equivalent forms will be called
equivalent. It is clear that equivalent surfaces are isomorphic.

\medskip

\noindent{\bf Main theorem.} {\em
Let $a_1,a_2,a_3$ be non-zero rational numbers. Let $X\subset\P^3_\Q$ be the surface
$$x_0^4+a_1x_1^4+a_2x_2^4+a_3x_3^4=0.$$
{\rm (1)} Let $k=\Q(i)$ and let $X_k=X\times_\Q k$. 
If $X_k$ is equivalent to
the diagonal quartic surface with coefficients $(a_1,a_2,a_3)=(1,2,-2)$, then
the $2$-primary torsion subgroup of $\Br(X_k)/\Br_1(X_k)$ is $\Z/2$; otherwise this subgroup is $0$.

\smallskip

\noindent{\rm (2)} 
If $X$ is equivalent to the diagonal quartic surface with coefficients $(1,2,-2)$ or $(1,8,-8)$, then
the $2$-primary torsion subgroup of $\Br(X)/\Br_1(X)$ is $\Z/2$; otherwise this subgroup is~$0$.}

\medskip

The two exceptional surfaces in (2) are 
non-equivalent over $\Q$ but equivalent, hence isomorphic over $\Q(i)$.
We do not know whether they are isomorphic over~$\Q$. 
Both have rational points over $\Q$.
In Proposition \ref{comp} we show that multiplying two of the coefficients by
$-4$ gives rise to surfaces with isomorphic Galois representations
in the second \'etale cohomology group
with coefficients $\Q_2$ and $\Z_\ell$ for any odd prime $\ell$.
By analogy with the Tate module of abelian varieties 
one may conjecture that these surfaces are related by an
isogeny whose degree is a power of $2$.

By the main theorem, only finitely many diagonal quartic surfaces 
with coefficients in $\Q$, considered
over $\Q$ or $\Q(i)$, up to isomorphism, have non-zero
2-torsion in the transcendental Brauer group. In contrast, there are infinitely many
non-equivalent 
diagonal quartic surfaces over $\Q$ with odd order torsion in the transcendental Brauer group.
Indeed, by \cite[Thm.~1.1]{IS} the {\em odd} torsion subgroup of $\Br(X)/\Br_1(X)$ is $0$, unless
$-3a_1a_2a_3$ is a fourth power or $-4$ times a fourth power in $\Q$ when this subgroup is $\Z/3$, 
or $125a_1a_2a_3$ is a fourth power or $-4$ times a fourth power in $\Q$
when this subgroup is $\Z/5$. (When $\Q$ is replaced by $\Q(i)$ we have $(\Z/3)^2$ and
$(\Z/5)^2$ instead of $\Z/3$ and $\Z/5$, respectively.)
Thus we now have a complete classification of the transcendental
Brauer groups of diagonal quartic surfaces with coefficients in $\Q$ over $\Q$ and $\Q(i)$. 
Note that the method of our paper works for all primes $\ell$, so our Proposition \ref{Gal}
directly implies
\cite[Thm.~1.1]{IS}, which was proved by a different method that works only for odd $\ell$.

As a consequence, we show non-existence of a 2-primary Brauer--Manin obstruction to 
the Hasse principle alluded to in \cite[p. 660]{ISZ}.

\medskip

\noindent{\bf Corollary} {\em Let $X$ be a diagonal quartic surface over $\Q$
with coefficients $a_1,a_2,a_3$ such that
$a_1a_2a_3$ is a square in $\Q$. Then $\Br(X)=\Br_1(X)$.}

\medskip

We can be even more precise: 
carrying out calculations of step (c) we describe the structure of the extension
\begin{equation}
0\lra \Br_1(X)/\Br_0(X)\lra \Br(X)/\Br_0(X) \lra \Br(X)/\Br_1(X)\lra 0. \label{ext}
\end{equation}

\medskip

\noindent{\bf Supplement to the main theorem} (1) {\em Let $X_k$ be the diagonal quartic surface with coefficients $(1,2,-2)$
over $k=\Q(i)$. Then the exact sequence $(\ref{ext})$ for $X_k$ over $k$ is the extension
$$0\lra \Z/2\times\Z/4\lra (\Z/4)^2\lra \Z/2\lra 0.$$
\noindent{\rm (2)} Let $X$ be the diagonal quartic surface with coefficients $(1,2,-2)$ or $(1,8,-8)$
over $\Q$. Then the exact sequence $(\ref{ext})$ for $X$ over $\Q$
is the extension}
$$0\lra\Z/4\lra \Z/8\lra \Z/2\lra 0.$$

Combined with Bright's classification of algebraic Brauer groups of diagonal quartic surfaces \cite{Bright},
this completes the classification of the Brauer groups of diagonal quartic surfaces 
with coefficients in $\Q$ over the ground fields $\Q$ and $\Q(i)$.

The supplement illustrates a difficulty in lifting Galois-invariant elements of $\Br(\ov X)[2]$ to $\Br(X)$
already apparent in work of Ieronymou \cite{I}:
an element in the image of $\Br(X)\to \Br(\ov X)$ may not lift
to an element of $\Br(X)$ of the same order. This difficulty does not arise for
odd order elements because $\Br_1(X)/\Br_0(X)$ is always a 2-group, as $\Ga$ acts on $\Pic(\ov X)$
via a finite 2-group.

The 2-primary torsion subgroup of $\Br(\ov X)$ contains a unique
non-zero element invariant under the group of diagonal automorphisms of $\ov X$. 
This element is Galois-invariant (Proposition \ref{nonzero}, Corollary \ref{coco}).
Our calculations in step (b) show that the map $\Br(X)\to\Br(\ov X)$ is almost always zero on
the 2-primary torsion subgroups, but when $\Br(X)/\Br_1(X)$ is non-zero in the exceptional cases above,
it is generated by this element. This answers Question 1 raised at the end of \cite{Sko17}.

The first named author has successfully applied the methods developed in this paper
to diagonal surfaces of degrees $5$ and $6$.
These methods also work for diagonal surfaces in weighted projective spaces,
see \cite{GLN, G}.

Let us outline the structure of the paper. In the preliminary Section 1
we develop necessary cohomological machinery and prove a finiteness result
for the Brauer group of diagonal surfaces over fields finitely generated over $\Q$.
In Section 2 we recall the description of the integral primitive cohomology of the complex Fermat
surface of degree $d$, including the cup-product
and Hodge decomposition, and explain how to recover
the full cohomology lattice from the primitive cohomology
sublattice. We also recall a known description of the associated
representation of the absolute Galois group of the cyclotomic field $\Q(\mu_d)$.
In Section 3 we specialise to the case $d=4$ and compute the transcendental lattice
of the quartic Fermat surface $F$ and the representation of the absolute Galois group of $\Q$
in the $\ell$-adic Tate module of $\Br(\ov F)$.
This allows us to carry out step (a). Although our general methods do not rely on the assumption that
the geometric Picard group of the Fermat surface is integrally generated by lines 
(which, by a theorem of Degtyarev \cite{D}, holds if and only if $d\leq 4$ or $d$ is coprime to 6;
for the case $d=4$ see \cite{PSD}),
our calculations in step (b) become simpler if we use this fact. We prove the main theorem
in Section~4 (see Theorems \ref{main1} and \ref{main2}); 
details of computations are collected in Appendix~A. The computational
proof of the supplement can be found in Appendix B.

This paper owns its existence to many happy discussions of the second named author
with Peter Swinnerton-Dyer, and to 
Peter's lifelong interest in diagonal quartics. We dedicate this paper to his memory.

We are grateful to Martin Bright, Rachel Newton, Domenico Valloni, Olivier Wittenberg and Yuri Zarhin 
for their questions and useful discussions.
The integral basis of lines on the Fermat quartic surface used in the appendix was 
kindly provided to us by Peter Swinnerton-Dyer.

The final part of the work on this paper was done at the Institut Henri Poincar\'e in Paris.
We thank the IHP for hospitality and support.

\section{Surfaces with torsion-free Picard group}

\subsection{Cohomological tools} \label{cohtools}

We start by recalling some basic homological algebra from \cite[\S 3.1]{CTS}.
Let $\mathcal A$, $\mathcal B$ be abelian categories
and let $G\colon\mathcal A\to\mathcal B$ be an additive functor. Let
\begin{equation}
0\lra A\lra B\lra C \lra 0 \label{ABC}
\end{equation}
be an exact sequence in $\mathcal A$.
Using injective resolutions one constructs a 2-term complex 
$\mathcal E=[E^0\to E^1]$ in $\mathcal B$ together with maps of complexes
\begin{equation}
\left(\tau_{[1,2]}G(B^\bullet)\right)[1]\longleftarrow {}\mathcal E\lra [(R^1G)C\to (R^2G)A],
\label{x1}
\end{equation}
where $B^\bullet$ is an injective resolution of $B$ and $\tau_{[1,2]}$ is the natural truncation. (See the diagram in \cite{CTS}, the proof of Lemma 3.2.) It satisfies the following properties.

\smallskip

(1) The leftward arrow induces an isomorphism of cohomology groups 
in degree 0. In degree 1 it induces the natural map from 
${\rm Coker}[(R^1G)C\to (R^2G)A]$ to $(R^2G)B$.

(2) The rightward arrow induces an isomorphism in degree 1. In degree 0 it induces
the natural map from $(R^1G)B$ to $\Ker[(R^1G)C\to (R^2G)A]$.

Let $k$ be a field with separable closure $\bar k$ and Galois group $\Ga=\Gal(\bar k/k)$. 
Let $p\colon X\to\Spec(k)$ be a variety over $k$. We write $\ov X=X\times_k\bar k$. 
The structure morphism $p$ defines the derived functor ${\bf R}p_*$
from the derived category $\mathcal D(X)$ of bounded below complexes of \'etale sheaves 
of abelian groups on $X$ to the derived category $\mathcal D(k)$
of bounded below complexes of continuous discrete $\Ga$-modules.

\ble \label{paris3}
Let $p\colon X\to\Spec(k)$ be a proper, geometrically reduced and geometrically connected 
variety over a field $k$.
There is a natural map $$\Br(X)/\Br_0(X)\lra \mathbb H^2(k, \tau_{[1,2]}{\bf R}p_*\G_{m,X}),$$
which is an isomorphism when $k$ is a number field.
\ele
{\em Proof.} The assumption on $X$ implies $\H^0(\ov X,\O_X)\cong\bar k$, hence
$\H^0(\ov X,\G_m)\cong \bar k^*$.
The derived functor of $\H^0(k,\cdot)$ is a functor from $\mathcal D(k)$
to the derived category of bounded below complexes of abelian groups. Applying it to the exact triangle 
$$\tau_{[0]}{\bf R}p_*\G_{m,X}\lra \tau_{[0,2]}{\bf R}p_*\G_{m,X}\lra
\tau_{[1,2]}{\bf R}p_*\G_{m,X}$$
gives an exact sequence
$$\Br(k)\lra\Br(X)\lra \mathbb H^2(k, \tau_{[1,2]}{\bf R}p_*\G_{m,X})\lra \H^3(k,\bar k^\times).$$
If $k$ is a number field, then $\H^3(k,\bar k^\times)=0$. $\Box$

\medskip

We now assume that $k$ has characteristic zero and that $X$ is a smooth, proper and geometrically
integral variety over $k$ such that $\Pic(\ov X)$ is torsion-free.
Let (\ref{ABC}) be the Kummer exact sequence
$$1\lra \mu_n\lra \G_m\lra \G_m \lra 1,$$
where $n\geq 1$. This sequence shows that
$\H^1(\ov X,\mu_n)=0$ for any $n\geq 1$,
thus in this case the rightward arrow in (\ref{x1}) is a quasi-isomorphism.
In the derived category $\mathcal D(k)$ we obtain a morphism 
$$
[\Pic(\ov X)\to \H^2(\ov X,\mu_n)] \lra (\tau_{[1,2]}{\bf R}p_*\G_{m,X})[1]
$$
inducing the isomorphism $n\Pic(\ov X)\tilde\lra\Pic(\ov X)$, $nx\mapsto x$, 
in degree 0 and the natural inclusion $\Br(\ov X)[n]\to \Br(\ov X)$ in degree 1.
If $\ell$ is a prime and $m$ is a positive integer, then setting $n=\ell^m$ and passing to the inductive limit
as $m\to\infty$ we obtain a morphism
$$
[\Pic(\ov X)\otimes\Z[\ell^{-1}]\to \H^2(\ov X,\Q_\ell/\Z_\ell(1))] \lra (\tau_{[1,2]}{\bf R}p_*\G_{m,X})[1].
$$
In degree 0 it induces the identity on $\Pic(\ov X)$ and in degree 1 the natural inclusion
$\Br(\ov X)\{\ell\}\hookrightarrow\Br(\ov X)$, where $\Br(\ov X)\{\ell\}$ is the $\ell$-primary torsion
subgroup, see \cite[II, Thm.~3.1]{GB}. Since $X$ is smooth, $\Br(\ov X)$ is a torsion group
\cite[II, Prop.~1.4]{GB}. Thus
summing over all primes $\ell$ we obtain a quasi-isomorphism
\begin{equation}
[\Pic(\ov X)\otimes\Q\to \H^2(\ov X,\Q/\Z(1))] \tilde\lra (\tau_{[1,2]}{\bf R}p_*\G_{m,X})[1].
\label{paris2}\end{equation}

Let $X$ be a smooth, projective, geometrically integral {\em surface} over a field $k\subset\C$
such that $\Pic(X_\C)$ is torsion-free. Then $\H^2(X_\C,\Z(1))$ is a torsion-free abelian group
and the first Chern class identifies $\Pic(X_\C)$ with a saturated subgroup of $\H^2(X_\C,\Z(1))$ 
(that is, the quotient is torsion-free).
The cup-product defines a unimodular symmetric bilinear pairing
$$\H^2(X_\C,\Z(1)) \times \H^2(X_\C,\Z(1)) \lra \H^4(X_\C,\Z(2))=\Z.$$
Since $\H^2(X_\C,\Z(1))$ is torsion-free, the cup-product induces an isomorphism
$$\H^2(X_\C,\Z(1))\tilde\lra\Hom(\H^2(X_\C,\Z(1)),\Z).$$
The Hodge index theorem implies that the restriction of the cup-product to $\Pic(X_\C)$,
which coincides with the intersection pairing on $\Pic(X_\C)$, has non-zero discriminant.
Define the {\em transcendental lattice} $T(X_\C)$ as the orthogonal complement 
to $\Pic(X_\C)$ in $\H^2(X_\C,\Z(1))$ with respect to the cup-product. 
Thus $\Pic(X_\C)\cap T(X_\C)=0$, and $\Pic(X_\C)\oplus T(X_\C)$ has finite index in $\H^2(X_\C,\Z(1))$.
Since $\Pic(X_\C)$ is a saturated subgroup of $\H^2(X_\C,\Z(1))$, it is the orthogonal complement to
$T(X_\C)$. Using the unimodularity of the cup-product on $\H^2(X_\C,\Z(1))$
we obtain canonical isomorphisms of finite abelian groups
$$\frac{\Hom(\Pic(X_\C),\Z)}{\Pic(X_\C)} \cong \frac{\H^2(X_\C,\Z(1))}{\Pic(X_\C)\oplus T(X_\C)} \cong
\frac{\Hom(T(X_\C),\Z)}{T(X_\C)}.$$
Let us call this group $\Delta_X$, or simply $\Delta$, if there is no confusion. 
We write $\Delta=\oplus_\ell\Delta_\ell$, where the order of $\Delta_\ell$
is a power of a prime number $\ell$.
We deduce an exact sequence of abelian groups
$$0\lra\Delta_\ell\lra \big(\Pic(X_\C)\otimes\Q_\ell/\Z_\ell\big)\oplus \big(T(X_\C)\otimes\Q_\ell/\Z_\ell\big) \lra
\H^2(X_\C,\Z(1))\otimes\Q_\ell/\Z_\ell\lra 0.$$
The group of connected components of the Picard scheme of $X\times_k K$,
where $K$ is an algebraically closed field extension of $k$,
does not depend on the choice of $K$. Thus, under the assumption that $\Pic(X_\C)$
is torsion-free, the groups $\Pic(\ov X)$ and $\Pic(X_\C)$
are canonically isomorphic. By the comparison theorem between classical and \'etale cohomology
we have an isomorphism $\H^2(X_\C,\Z(1))\otimes\Z_\ell\cong\H^2(\ov X,\Z_\ell(1))$,
compatible with the cycle class map and the cup-product. 
Thus $T(X_\C)\otimes\Z_\ell\cong T(\ov X)_\ell$,
where $T(\ov X)_\ell\subset \H^2(\ov X,\Z_\ell(1))$ is the orthogonal complement
to the (injective) image of $\Pic(\ov X)$ under the class map, and
the previous exact sequence
is canonically identified with the exact sequence of $\Ga$-modules
$$0\lra\Delta_\ell\lra \big(\Pic(\ov X)\otimes\Q_\ell/\Z_\ell\big)\oplus 
\big(T(\ov X)_\ell\otimes\Q_\ell/\Z_\ell\big) \lra
\H^2(\ov X,\Z_\ell(1))\otimes\Q_\ell/\Z_\ell\lra 0.$$
This gives rise to the commutative diagram of $\Ga$-modules
$$\xymatrix{\Delta_\ell \ar@{^{(}->}[r]\ar@{^{(}->}[d]&
T(\ov X)_\ell\otimes\Q_\ell/\Z_\ell\ar@{^{(}->}[d]^{[-1]}\\
\Pic(\ov X)\otimes\Q_\ell/\Z_\ell \ar@{^{(}->}[r] &\H^2(\ov X,\Z_\ell(1))\otimes\Q_\ell/\Z_\ell
}$$
Summing over all primes $\ell$ we obtain a commutative diagram of $\Ga$-modules
$$\xymatrix{\Delta\ar@{^{(}->}[r]\ar@{^{(}->}[d]&T(X_\C)\otimes\Q/\Z\ar@{^{(}->}[d]^{[-1]}\\
\Pic(X_\C)\otimes\Q/\Z\ar@{^{(}->}[r] &\H^2(X_\C,\Z(1))\otimes\Q/\Z
}$$
Using the non-degeneracy of the intersection pairing on $\Pic(\ov X)$ we identify
the $\Ga$-modules $\Pic(\ov X)\otimes\Q$ and $\Hom(\Pic(\ov X),\Q)$.
This gives rise to the commutative diagram of $\Ga$-modules
$$\xymatrix{\Hom(\Pic(X_\C),\Z)\ar@{->>}[r]\ar[d]&\Delta\ar[d]\\
\Pic(X_\C)\otimes\Q\ar@{->>}[r]&\Pic(X_\C)\otimes\Q/\Z
}$$
Finally, we compose horizontal arrows in the two last diagrams to obtain a commutative diagram of $\Ga$-modules
\begin{equation}
\begin{split}
\xymatrix{\Hom(\Pic(X_\C),\Z)\ar[r]\ar[d]&T(X_\C)\otimes\Q/\Z\ar[d]^{[-1]}\\
\Pic(X_\C)\otimes\Q\ar[r] &\H^2(X_\C,\Z(1))\otimes\Q/\Z
} \label{dia}
\end{split}
\end{equation}

Since $X$ is a surface, by Poincar\'e duality $\H^2(\ov X,\Z_\ell)_{\rm tors}=0$ 
implies $\H^3(\ov X,\Z_\ell)_{\rm tors}=0$, for all $\ell$. 
By the Kummer sequence this implies that $\Br(\ov X)\cong\Br(X_\C)$ is divisible
and the $\ell$-primary torsion subgroup $\Br(\ov X)\{\ell\}$ is the quotient of 
$\H^2(\ov X,\Z_\ell(1))\otimes\Q_\ell/\Z_\ell$ by $\Pic(\ov X)\otimes\Q_\ell/\Z_\ell$,
see \cite{CTS}, diagram (26), p.~157.
Taking the direct sum over all primes $\ell$
we obtain a canonical isomorphism of $\Ga$-modules
\begin{equation}
\Hom_\Z(T(X_\C),\Q/\Z) \tilde\lra \Br(\ov X). \label{e1}
\end{equation}
For future use we point out that (\ref{e1}) gives rise to an isomorphism of $\Ga$-modules
\begin{equation}
\Hom_\Z(T(X_\C),\Z/\ell^m) \tilde\lra \Br(\ov X)[\ell^m], \label{e19}
\end{equation}
for any prime $\ell$ and any positive integer $m$.

\bpr \label{le2}
Let $X$ be a smooth, projective, geometrically integral surface over a subfield $k\subset\C$
such that $\Pic(X_\C)$ is torsion-free. Then the map of
horizontal $2$-term complexes of $\Ga$-modules 
in $(\ref{dia})$ is a quasi-isomorphism. 
Each of these complexes concentrated in degrees $0$ and $1$
represents $(\tau_{[1,2]}{\bf R}p_*\G_{m,X})[1]$ in $\mathcal D(k)$.
\epr
{\em Proof.} The first statement follows from the construction of (\ref{dia}). Indeed,
the induced map in degree 0 is the identity on $\Pic(X_\C)=\Pic(\ov X)$. 
The induced map in degree 1 is the canonical isomorphism (\ref{e1}).

Since $\H^3(X_\C,\Z)_{\rm tors}=0$, the exact sequence
$$0\lra\Z\lra\Q\lra\Q/\Z\lra 0$$
gives a natural isomorphism 
$$\H^2(X_\C,\Z(1))\otimes\Q/\Z \tilde\lra \H^2(X_\C,\Q/\Z(1)).$$
Now the quasi-isomorphism (\ref{paris2}) gives that the bottom complex of (\ref{dia})
is quasi-isomorphic to $(\tau_{[1,2]}{\bf R}p_*\G_{m,X})[1]$. $\Box$

\brem{\rm For completeness we note that 
$\tau_{[0,1]}{\bf R}p_*\G_{m,X}$ is represented by
the complex of $\Ga$-modules $\div\colon\bar k(X)^\times\to \Div(\ov X)$, 
where $\div$ sends a function to its divisor, see \cite[Lemma 2.3]{BvH09}.}
\erem

\bco \label{1.4}
Let $X$ be a smooth, projective, geometrically integral surface over a 
number field $k$
such that $\Pic(\ov X)$ is torsion-free. Then the group $\Br(X)/\Br_0(X)$ is canonically isomorphic
to the Galois hypercohomology group
$$\mathbb H^1(k,\left[\Hom(\Pic(\ov X),\Z)\to T(X_\C)\otimes\Q/\Z\right]).$$
\eco
{\em Proof.} This follows from Proposition \ref{le2} in view of Lemma \ref{paris3}. $\Box$

\medskip

Recall that we have exact sequences of $\Ga$-modules
\begin{equation}
0\lra\Delta\lra T(X_\C)\otimes_\Z\Q/\Z\lra \Br(\ov X)\lra 0\label{three}
\end{equation}
and
\begin{equation}
0\lra\Pic(\ov X)\lra \Hom_\Z(\Pic(\ov X),\Z)\lra \Delta\lra 0.\label{one}
\end{equation}

\bco \label{1}
Let $X$ be a smooth, projective, geometrically integral surface over a number field $k$
such that $\Pic(\ov X)$ is torsion-free. Then
the image of the natural map $\Br(X)\to \Br(\ov X)$
is equal to the kernel of the composition 
$$\Br(\ov X)^{\Ga}\stackrel{\partial_1}\lra\H^1(k,\Delta)
\stackrel{\partial_2}\lra\H^2(k,\Pic(\ov X)),$$
where $\partial_1$ is defined by $(\ref{three})$ and $\partial_2$
is defined by $(\ref{one})$.
\eco
{\em Proof.} This follows from Proposition \ref{le2}
using that $\H^3(k,\bar k^\times)=0$ for a number field $k$. See also \cite[Prop.~4.1]{CTS}. $\Box$

\medskip

\bco \label{1.5}
Let $X$ be a smooth, projective, geometrically integral surface over a number field $k$
such that $\Pic(\ov X)$ is torsion-free. Suppose that the image of 
the natural map $\Br(X)\to\Br(\ov X)^\Ga$ is contained in $\Br(\ov X)[n]$
for some $n\geq 1$. Let $T_n\subset T(X_\C)\otimes\Q/\Z$ be the inverse image of $\Br(\ov X)[n]$
in $(\ref{three})$. 
Let $\Ga_K$ be the kernel of the natural map 
$$\Ga\lra \Aut(\Pic(\ov X))\times\Aut(T_n),$$
and let $K=(\bar k)^{\Ga_K}$, $G=\Gal(K/k)$. Then we have the following statements.

{\rm (a)} The image of the natural map $\Br(X)\to \Br(\ov X)^\Ga$
is equal to the kernel of the composition
$$\Br(\ov X)[n]^{G}\stackrel{\partial_1}\lra\H^1(G,\Delta)
\stackrel{\partial_2}\lra\H^2(G,\Pic(\ov X)).$$

{\rm (b)} 
The group $\Br(X)/\Br_0(X)$ is canonically isomorphic
to the hypercohomology group
\begin{equation}
\mathbb H^1(G,\left[\Hom(\Pic(X_\C),\Z)\to T_n\right]). \label{hyper}
\end{equation}
\eco
{\em Proof.} (a) 
It is clear that the action of $\Ga$ on the terms of (\ref{one}) and on the terms of
$$0\lra\Delta\lra T_n\lra \Br(\ov X)[n]\lra 0$$
factors through $G$. 
We have $\Pic(\ov X)^{\Ga_K}=\Pic(\ov X)$, which is a finitely generated
free abelian group, hence $\H^1(K,\Pic(\ov X))=0$. 
From the Hochschild--Serre
spectral sequence we deduce the injectivity of the inflation map
$$\H^2(G,\Pic(\ov X))\lra \H^2(k,\Pic(\ov X)).$$ Now (a) follows from Corollary \ref{1}.

(b) Let $\K$ be the complex $\Hom(\Pic(\ov X),\Z)\to T(X_\C)\otimes\Q/\Z$
and let $\K_n$ be the complex $\Hom(\Pic(X_\C),\Z)\to T_n$.
Using Proposition \ref{le2} we obtain a commutative diagram, whose lower part is constructed using inflation maps:
$$\xymatrix{0\ar[r]&\H^1(k,\Pic(\ov X))\ar[r]&\mathbb{H}^1(k,\K)\ar[r]&\Br(\ov X)^\Ga\ar[r]&
\H^2(k,\Pic(\ov X))\\
0\ar[r]&\H^1(k,\Pic(\ov X))\ar[r]\ar[u]^\cong&\mathbb{H}^1(k,\K_n)\ar[r]\ar[u]&
\Br(\ov X)[n]^\Ga\ar[r]\ar@{^{(}->}[u]&\H^2(k,\Pic(\ov X))\ar[u]^\cong\\
0\ar[r]&\H^1(G,\Pic(\ov X))\ar[r]\ar[u]^\cong&\mathbb{H}^1(G,\K_n)\ar[r]\ar[u]&
\Br(\ov X)[n]^G\ar[r]\ar[u]^\cong&\H^2(G,\Pic(\ov X))\ar@{^{(}->}[u]}$$
A diagram chase gives canonical isomorphisms 
$$\mathbb{H}^1(G,\K_n)\tilde\lra\mathbb{H}^1(k,\K_n)\tilde\lra\mathbb{H}^1(k,\K)\tilde\lra\Br(X)/\Br_0(X),$$
the last one given by Corollary \ref{1.4}. $\Box$

\subsection{Brauer group of varieties dominated by products of curves}

A smooth, projective and geometrically integral variety $X$ over a field $k$ is called
a variety {\em dominated by a product of curves} if there is a dominant rational map
from a product of geometrically integral $\bar k$-curves to $\ov X$.

\bpr \label{19nov}
Let $k$ be a field finitely generated over $\Q$. 
If $X$ is a smooth, projective and geometrically integral 
variety over $k$ dominated by a product of curves, then $\Br(\ov X)^\Ga$ is finite.
\epr
{\em Proof.} If $V$ and $W$ are smooth, projective and geometrically integral varieties
over a field $k$ which is finitely generated over $\Q$, then the cokernel of the natural
map $\Br(\ov V)^\Ga\oplus\Br(\ov W)^\Ga\to\Br(\ov V\times\ov W)^\Ga$
is finite by \cite[Thm.~A]{SZ14}. The Brauer group of a 
curve over an algebraically closed field is zero \cite[III, Cor.~1.2]{GB}, \cite[Thm.~5.6.1 (ii)]{CTS21}.
Thus if $Z$ is a product of smooth, projective and geometrically integral curves over $k$, 
then $\Br(\ov Z)^\Ga$ is finite.

Replacing $k$ by a finite field extension we obtain
geometrically integral $k$-curves $C_1,\ldots,C_n$ such that there is a
dominant rational map $f$ from $S=\prod_{i=1}^n C_i$ to $X$ defined over $k$. Let $d=\dim(X)$.
There is a dense open subset $U\subset S$ such that the restriction of $f$ to $U$
is a smooth morphism $U\to X$. After another finite extension of $k$ we can find a $k$-point
$P\in U(k)$. The induced map of tangent spaces
$f_*\colon T_{S,P}\to T_{X,f(P)}$ is surjective, so we can choose
a $d$-element subset $I\subset\{1,\ldots,n\}$ such that $f_*(T_{S_I,P})=T_{X,f(P)}$,
where $S_I\cong\prod_{i\in I}C_i$ is the $k$-fibre of the projection
$S\to\prod_{i\notin I}C_i$ which contains $P$.
Thus $f$ restricts to a dominant, generically finite rational map from a product of $d$ geometrically integral
$k$-curves to $X$.

Since ${\rm char}(k)=0$, we can assume that there is a smooth, projective and geometrically integral
variety $Y$ over $k$, a birational morphism
$Y\to Z$, where $Z$ is a product of smooth, projective and geometrically integral curves,
and a dominant, generically finite morphism $f\colon Y\to X$. By the birational invariance of the Brauer group
\cite[III, \S 7]{GB}, \cite[Cor.~6.2.11]{CTS21}
we have $\Br(\ov Y)^\Ga\cong\Br(\ov Z)^\Ga$. A general theorem of Grothendieck says that the natural
map $\Br(\ov X)\hookrightarrow\Br(\bar k(X))$ is injective, see \cite[III, Cor.~6.2]{GB},
\cite[Thm.~3.5.5]{CTS21}.
The standard restriction-corestriction 
argument then gives that the kernel of $f^*\colon\Br(\ov X)\to\Br(\ov Y)$
is killed by the degree $[\bar k(Y):\bar k(X)]$. Since $\Br(\ov X)$ is a torsion group of
cofinite type \cite[II, Cor.~3.4]{GB}, \cite[Prop.~5.2.9]{CTS21}, this kernel is finite. 
Hence $\Br(\ov X)^\Ga$ is finite. $\Box$

\medskip

\bco
Let $k$ be a field finitely generated over $\Q$.
Let $f(t)$ and $g(t)$ be separable polynomials of degree $d\geq 2$.
Let $F(x,y)$ and $G(x,y)$ be homogeneous forms of degree $d$ such that $f(t)=F(t,1)$ and $g(t)=G(t,1)$.
Let $X\subset \P^3_k$ be the surface with equation $F(x,y)=G(z,w)$, for example, a diagonal surface.
Then the Brauer group $\Br(X)$ is finite modulo $\Br(k)$.
\eco
{\em Proof.} Consider the smooth plane curves $C$ and $D$
given by $z^d=F(x,y)$ and $u^d=G(v,w)$, respectively.
Passing to the quotient by the action of $\mu_d$ on $C\times D$, which multiplies $z$ and $u$ by
the same root of unity, gives a dominant rational map from $C\times D$ to $X$. Thus $X$ is a DPC surface.
Since $\Pic(\ov X)$ is torsion-free \cite[XI.1.8]{SGA},
the group $\Br_1(X)/\Br_0(X)$ is finite. The finiteness of $\Br(\ov X)^\Ga$ follows 
from Proposition \ref{19nov}. $\Box$

\section{Cohomology of Fermat surfaces}

\subsection{Topological background}

Let $X$ be a smooth irreducible subvariety of $\P^{n+1}_\C$ of dimension $n$.
Let $Z\subset X$ be a smooth hyperplane section and let $U=X\setminus Z$.
The long exact sequence of relative homology groups and the long exact sequence of
cohomology groups with compact 
support form the following commutative diagram, where vertical arrows
are the Poincar\'e duality isomorphisms (see, e.g., the diagram in the proof of \cite[Prop.~3.46]{Hatch}):
\begin{equation}
\begin{split}\xymatrix{
\H_{n+1}(X,U,\Z)\ar[r]& \H_n(U,\Z)\ar[r]& \H_n(X,\Z)\ar[r]& \H_n(X,U,\Z)\\
\H^{n-1}(Z,\Z)\ar[r]\ar[u]^\cong& \H^n_c(U,\Z)\ar[r]\ar[u]^\cong
& \H^n(X,\Z)\ar[r]\ar[u]^\cong& \H^n(Z,\Z)\ar[u]^\cong
}\label{may1}
\end{split}
\end{equation}
Next, we have a commutative diagram, where the maps $h$ are defined by evaluating
cochains on chains and the isomorphisms are dual to the Poincar\'e duality isomorphisms:
$$\xymatrix{
\H^n(U,\Z)\ar[r]^{h\quad\quad}&\Hom(\H_n(U,\Z),\Z)\ar[r]^{\sim\ }&
\Hom(\H^n_c(U,\Z),\Z)\\
\H^n(X,\Z)\ar[r]^{h\quad\quad}\ar[u]&\Hom(\H_n(X,\Z),\Z)\ar[r]^{\sim\ }\ar[u]
&\Hom(\H^n(X,\Z),\Z)\ar[u]
}$$
The commutativity of this diagram implies compatibility of the cup-product pairings
$$\begin{array}{ccccc}
\H^n(U,\Z)&\times &\H^n_c(U,\Z)&\lra&\Z\\
\uparrow&&\downarrow&&||\\
\H^n(X,\Z)&\times &\H^n(X,\Z)&\lra&\Z
\end{array}$$
Modulo torsion subgroups, these pairings are perfect dualities, cf.~\cite[Prop.~3.38]{Hatch}.
The natural map $\H^n_c(U,\Z)\to \H^n(U,\Z)$ gives rise to a bilinear form on
$\H^n_c(U,\Z)$.  The Poincar\'e duality isomorphism
$\H^n_c(U,\Z)\cong \H_n(U,\Z)$ identifies this form with the intersection products
of homology classes, see \cite[pp.~43--60]{GH}.

\subsection{Primitive cohomology and the cup-product} \label{prim_coh}

Let $F\subset \P^3_\Q$ be the Fermat surface $$x_0^d+x_1^d+x_2^d+x_3^d=0.$$
The abelian group $H=\H^2(F_\C,\Z(1))$ is finitely generated and free, so
the cup-product defines a unimodular symmetric bilinear form $\langle x,y\rangle$ on $H$.
Let $L\in H$ be the hyperplane section class. 
Let $P=P^2(F_\C,\Z(1))\subset H$ be the orthogonal complement to $L$
with respect to the cup-product pairing. 
Thus $P$ is the kernel of the map $H\to\Z$ sending $x\in H$ to $\langle x, L\rangle$, so
$P$ is a saturated subgroup of $H$. Since $\langle L,L\rangle=d>0$, the restriction of the cup-product pairing
$H\times H \to\Z$ to $P$ is non-degenerate.

Let $\Lambda\in H$ be the class of a line 
contained in $F_\C$. Since $\langle L, \Lambda\rangle=1$, we see that $L$ is a primitive element of $H$.
(This holds for arbitrary smooth complete intersection surfaces in a projective space, see \cite[XI.1.8]{SGA}.)
Our map $H\to\Z$ sends $\Lambda$ to 1, so we have a direct sum decomposition
of abelian groups $H= \Z\Lambda\oplus P$.

Let $G$ be the group $\Q$-scheme $(\mu_d)^3$. We write $G$ for $G(\ov\Q)$.
The choice of a primitive $d$-th root $\e$ of 1 defines 
an isomorphism of groups $G\simeq (\Z/d)^3$ so that
the elements of $G$ are monomials $u_1^a u_2^b u_3^c$, where 
$$u_1=(\e,1,1),\quad u_2=(1,\e,1),\quad u_3=(1,1,\e),$$
and $(a,b,c)\in (\Z/d)^3$.
The group $G$ acts on $F_\C$ by automorphisms so that $u_i$ multiplies $x_i$ by $\e$ 
and leaves the other coordinates unaltered.

Let $U\subset F$ be the complement to the smooth hyperplane section $Z\subset F$
given by $x_0=0$. We are in the situation of the previous section, with $X=F_\C$.
It is clear that the group $G\simeq(\Z/d)^3$ acts on all the terms of diagram (\ref{may1})
so that the rows are exact sequences of $G$-modules. The Poincar\'e duality isomorphisms
are defined by taking the cap-product with the fundamental class in $\H_4(F_\C,\Z)\cong\Z$,
given by the orientation of $F_\C$.
The fundamental class is fixed by the action of $G$,
hence the Poincar\'e duality isomorphisms are isomorphisms of $\Z[G]$-modules. 

From (\ref{may1}) we thus obtain an exact sequence of $G$-modules
\begin{equation}
0\lra\H_1(Z_\C,\Z)\lra \H_2(U_\C,\Z)\lra \H_2(F_\C,\Z)\lra \H_0(Z_\C,\Z)\cong\Z, \label{Gysin}
\end{equation}
where we used that $\H_1(F_\C,\Z)=0$. If we ignore the twists, then
the last arrow is the map $H\to\Z$ sending $x\in H$ to $\langle x, L\rangle$, so its kernel is $P$. 
Thus (\ref{Gysin}) gives rise to an exact sequence of $G$-modules
\begin{equation}
0\lra\H_1(Z_\C,\Z)\lra \H_2(U_\C,\Z)\lra P\lra 0, \label{Gysin_prim}
\end{equation}
where we ignore the twists.
(The twist of $P$ plays no role
in the description of the cup-product pairing and the action of $G$ on $P$, but will be essential
for the determination of the action of complex conjugation, see Remark \ref{compl}.)

Since the cup-product pairing on $P$ is non-degenerate, the kernel of the intersection pairing
on $\H_2(U_\C,\Z)$ is the image of $\H_1(Z_\C,\Z)$. 

\smallskip

Let us recall a description of the $G$-module $\H_2(U_\C,\Z)$ together with the intersection pairing on it,
obtained by Pham in \cite[Thm.~1]{Ph}. The affine Fermat surface $U\subset\A^3_k$ is given by
$$x_1^d+x_2^d+x_3^d=-1.$$
Let $S\subset\R^3$ be the standard 2-dimensional simplex 
$z_1+z_2+z_3=1$, $z_i\geq 0$, for $i=1,2,3$. Consider the map 
${\mathbf e}:S\to U_\C$ given by
$${\mathbf e}(z_1,z_2,z_3)=\left(\zeta_{2d}z_1^{1/d},\zeta_{2d}z_2^{1/d},\zeta_{2d}z_3^{1/d}\right),$$
where $\zeta_{2d}$ is a primitive
root of unity of degree $2d$ such that $\zeta_{2d}^2=\e$.
Then $$e=(1-u_1^{-1})(1-u_2^{-1})(1-u_3^{-1}){\mathbf e}\in\H_2(U_\C,\Z)$$ is a cycle.
Pham shows that the topological space of the simplicial complex generated by $g({\mathbf e})$,
for $g\in G$, is a $G$-equivariant deformation retract of $U_\C$. 
The resulting chain compex gives isomorphisms of $\Z[G]$-modules 
\begin{equation}
\H_2(U_\C,\Z)\cong\Z[G]e\cong\Z[u_1,u_1,u_3]/\big(\phi(u_1), \phi(u_2), \phi(u_3)\big),
\label{pham}
\end{equation}
where $\phi(x)$ is defined as follows:
\begin{equation}
\phi(x)= 1+x+x^2+\dots+x^{d-1}\ \in \ \Z[x], \quad \rho(x,y)=\sum_{0\leq i\leq j\leq  d-2}y^i x^j
\ \in \ \Z[x,y].\label{phi&rho}
\end{equation}
Consider the ring
$R_x:=\Z[x]/(\phi(x))\simeq\Z^{d-1}$. It is immediate to check that in $R_x$
we have $(1-x)\rho(1,x)=d$, so $d(1-x)^{-1}\in R_x$. In this notation, Pham's isomorphism (\ref{pham})
is $\H_2(U_\C,\Z)\cong R_{u_1}\otimes R_{u_2}\otimes R_{u_3}$.

\ble \label{loo}
We have an isomorphism of $\Z[G]$-modules
$$P\cong \Z[u_1,u_1,u_3]/I, \quad\text{where}\quad 
I=\left(\phi(u_1), \phi(u_2), \phi(u_3), \phi(u_1u_2u_3)\right).$$
\ele
{\em Proof.} 
In \cite[p.~337]{Ph} Pham defines a $\Z[G]$-valued form $(x|y)$ on $\H_2(U_\C,\Z)$,
as follows:
$$(x|y)=\sum_{g\in G}\langle x,gy\rangle g^{-1}.$$
This form is $\Z[G]$-linear in $y$ and $\Z[G]$-antilinear in $x$, with respect to the involution of 
the ring $\Z[G]$ that sends $g$ to $g^{-1}$.
Then he computes \cite[p.~340]{Ph} 
\begin{equation}
(e|e)=-(1-u_1)(1-u_2)(1-u_3)(1-(u_1u_2u_3)^{-1}). \label{ee}
\end{equation}
(Note that the right-hand side of (\ref{ee}) is invariant under the involution of $\Z[G]$.)
Since $x$ is invertible in $R_x$, and $1-x$ becomes invertible in $R_x\otimes\Q$, we deduce
that the kernel of the intersection product on $\H_2(U_\C,\Z)$
is the annihilator of $1-u_1u_2u_3$ in the ring $R_{u_1}\otimes R_{u_2}\otimes R_{u_3}$.
We need to prove that this ideal is generated by $\phi(u_1u_2u_3)$. 

This can be rephrased in terms
of group cohomology of the cyclic group $C_d$. Namely, we have an exact sequence of $C_d$-modules
$$ 0\lra \Z\lra \Z[C_d]\lra R_x\lra 0,$$ 
where $\Z$ is a trivial $C_d$-module. Multiplication by $u_1u_2u_3$ on 
$R_{u_1}\otimes R_{u_2}\otimes R_{u_3}$ makes it a $C_d$-module isomorphic to
$(\Z[C_d]/\Z)^{\otimes 3}$. The annihilator ideal of $1-u_1u_2u_3$ is the $C_d$-invariant
subgroup of $(\Z[C_d]/\Z)^{\otimes 3}$, whereas $\phi(u_1u_2u_3)$ is the norm element of
the group ring $\Z[C_d]$. Thus we need to show that the 0-th Tate cohomology group
$\widehat \H^0(C_d,(\Z[C_d]/\Z)^{\otimes 3})$ is trivial.

Taking the triple self-product of the exact complex $0\to\Z\to \Z[C_d]$
we obtain the following exact sequence of $C_d$-modules:
$$0\lra\Z\lra \Z[C_d]^{\oplus 3}\lra (\Z[C_d]^{\otimes 2})^{\oplus 3}\lra \Z[C_d]^{\otimes 3}
\lra (\Z[C_d]/\Z)^{\otimes 3}\lra 0.$$
The three terms in the middle are free $\Z[C_d]$-modules, hence all their Tate cohomology groups
are trivial. We obtain
$$\widehat \H^0(C_d,(\Z[C_d]/\Z)^{\otimes 3})\cong \H^3(C_d,\Z)\cong \Hom(C_d,\Z)=0.$$
This finishes the proof. $\Box$

\brem{\rm 
This result is stated in \cite[Cor. 2.2]{L}. The argument in the above proof
gives that $\H^0(C_d,(\Z[C_d]/\Z)^{\otimes(n+1)})$ is zero when $n$ is even,
and is isomorphic to $\Z/d$ when $n$ is odd. This second case is missed in \cite[Cor. 2.2]{L}
so the statement of {\em loc.~cit.} for odd $n$ needs to be modified.
}
\erem

\brem \label{compl}{\rm
The complex conjugation is an involution of the real manifold $F(\C)$.
Denote by $\tau$ the induced involutions of the homology and cohomology
groups of $F_\C$ and $U_\C$ with coefficients $\Z$. Then $\tau g=g^{-1}\tau$ for any $g\in G$.
From $\tau({\mathbf e})=(u_1 u_2 u_3)^{-1}{\mathbf e}$ we deduce that $\tau(e)=-e$.

Let $e'\in \H_2(F_\C,\Z)$ be the image of $e\in \H_2(U_\C,\Z)$.
Since the maps in (\ref{Gysin}) respect the action of $\tau$, we have $\tau(e')=-e'$.
The complex conjugation $\tau$ preserves the orientation of $F_\C$, because $F_\C$ has complex dimension 2,
hence preserves the fundamental class, so the Poincar\'e duality isomorphisms 
respect the action of~$\tau$. 

Now let us consider the twisted cohomology group $H=\H^2(F_\C,\Z(1))$. By definition, there is an isomorphism
$$\H^2(F_\C,\Z)\tilde\lra \H^2(F_\C,\Z(1))$$
sending $x$ to $2\pi i x$. Thus $H$ with its natural (twisted) action of $\tau$ is a subgroup of 
the complex vector space $\H^2(F_\C,\C)$ where the complex conjugation is the product of 
$\tau$, which is induced by the conjugation acting on $F(\C)$,
and the complex conjugation in $\C$. Let us denote by $\tau(1)$ this twisted action on $H$.
In particular, for $2\pi i e'\in H$ we have
\begin{equation}
\tau(1)(2\pi i e')=2\pi i e'.\label{conj}
\end{equation}
As $2\pi i e'$ generates the $\Z[G]$-module $P$, the isomorphism of Lemma \ref{loo} 
identifies $2\pi i e'\in P$ with the image of 
$1 \in \Z[G]$, so this isomorphism is compatible with the action of $\tau(1)$ on $\Z[G]$ which fixes 1 and sends
each $g\in G$ to $g^{-1}$.
(According to \cite[Cor.~2.2]{L},
permutations of $\{1,2,3\}$ act on $e'$ by the sign character.)}
\erem

The cup-product symmetric bilinear form on $P$ has the following description.
Define $\varphi\colon\mu_d\to\Z$ by
$\varphi(1)=1$, $\varphi(\e)=-1$, and $\varphi(\e^i)=0$ for $i=2,\ldots,d-1$.

\bpr \label{o1}
Under the isomorphism of $G$-modules $P\cong \Z[G]/I$,
the cup-product is the unique symmetric $G$-invariant bilinear form 
$\langle x,y\rangle\colon P\times P\to\Z$ such that
$$\langle u_1^a u_2^b u_3^c,1\rangle
=-\varphi(\e^{a})\varphi(\e^{b})\varphi(\e^{c})-
\varphi(\e^{-a})\varphi(\e^{-b})\varphi(\e^{-c}).$$
For any $r\in \Z[G]$ we have $\langle rx,y\rangle=\langle x,\bar r y\rangle$,
where $r\mapsto \bar r$ is the automorphism of the abelian group $\Z[G]$ that sends
$g$ to $g^{-1}$ for any $g\in G$.
\epr
{\em Proof.} This is obtain by a straightforward calculation from (\ref{ee}). 
The last statement follows from $G$-invariance of the form.
$\Box$

\medskip

We write $E=\Q(\mu_d)$. It is well known that $\Gal(E/\Q)\cong(\Z/d)^\times$.

Let $\widehat G=\Hom(G,\C^\times)\cong(\Z/d)^3$ be the group of characters. 
Write $\chi(u_1^a u_2^b u_3^c)=\e^{la+mb+nc}$, where $l, m, n$ are elements of $\{0,\ldots,d-1\}$.
Then 
$$E[G]=E[u_1,u_2,u_3]/(u_i^d-1)= \bigoplus_{\chi\in \widehat G}V_\chi$$
is a direct sum of 1-dimensional eigenspaces $V_\chi$ generated by the orthogonal idempotents
$$\alpha_\chi=\alpha_l(u_1)\alpha_m(u_2)\alpha_n(u_3)\in E[G], \quad\text{where}\quad\alpha_i(u)=
\frac{1}{d}\sum_{j=0}^{d-1}\e^{-ij}u^j.$$
Indeed, one readily checks that $\alpha_\chi$ is an eigenvector with eigenvalue $\chi$:
\begin{equation}
[g]\alpha_\chi=\chi(g)\alpha_\chi\quad \text{for any} \ \ g\in G, \label{29dec1}
\end{equation}
where we write $[g]$ for the natural action of $G$ on $E[G]$. 
Hence $\alpha_\chi\alpha_\phi$ is a multiple of $\delta_{\chi,\phi}$, where $\delta$ is the Kronecker delta.
Our normalisation is such that $\alpha_\chi^2=\alpha_\chi$.

\bco \label{co}
For any $\chi,\phi\in\widehat G$ we have
$\langle \alpha_\chi, \alpha_\phi\rangle=0$ if  $\phi\not=\chi^{-1}$, and
\begin{equation}\langle \alpha_\chi, \alpha_{\chi^{-1}}\rangle=
\langle \alpha_\chi^2, 1\rangle=\langle \alpha_\chi, 1\rangle=
-\frac{2}{d^3}\,{\rm Re}\,\left( (1-\e^l)(1-\e^m)(1-\e^n)\right).
\label{cup}
\end{equation}
\eco
{\em Proof.} This follows directly from Proposition \ref{o1}. $\Box$

\medskip

We note that
$I_E=I\otimes E$ is the direct sum of the spaces $V_\chi$
such that $\chi$ restricts trivially to any of the three coordinates $\mu_d\subset G$ 
(that is, $l=0$ or $m=0$ or $n=0$) or
to the diagonal $\mu_d\subset G$ (that is, $l+m+n$ is divisible by $d$). 
Let us denote by $S$ the complement to this set in $\widehat G$. In other words, $S$ is the set of triples
$(l,m,n)\in \{1,\ldots, d-1\}^3$ such that $l+m+n$ is not divisible by $d$.
Then 
$$P\otimes E=\bigoplus_{\chi\in S}V_\chi,$$
so that the rank of $P$ is $|S|$. For $d=2,3,4$ we have $|S|=1, 6, 21$, respectively.

The group $\Gal(E/\Q)\cong(\Z/d)^\times$ 
acts on $\widehat G=\Hom(G,E^\times)$ via its natural action on $E$,
which raises $\e$ to powers coprime to $d$. Note that this action preserves $S\subset\widehat G$.
For $\si\in\Gal(E/\Q)$ we have 
$\langle \alpha_\chi, \alpha_\phi\rangle=\langle \alpha_{\si(\chi)}, \alpha_{\si(\phi)}\rangle$.

\subsection{Hodge decomposition} \label{Hodge}

The Hodge decomposition $P\otimes \C=P^{1,-1}\oplus P^{0,0}\oplus P^{-1,1}$ can be described as follows.
Let $q$ be the function $S\to\{-1,0,1\}$ given by $q(\chi)=\left \lfloor{\frac{l+m+n}{d}}\right \rfloor -1$. 
Note that $q$ depends on the choice of $\e$.

\ble \label{o2}
The Hodge subspace $P^{p,q}\subset P\otimes \C=P^2(F_\C,\C(1))$ is the direct sum of $V_\chi\otimes_E\C$
for $\chi\in S$ such that $q(\chi)=q$.
\ele 
{\em Proof.} This is \cite[Thm.~I]{Sh79}. For a proof see \cite[Prop.~2.7]{ABB} based on \cite{L}.~$\Box$

\medskip

In particular, the $\C$-vector space $P^{p,q}$ comes from the $E$-vector space $\oplus_{q(\chi)=q}V_\chi$.

Let $S_\sharp\subset S$ be the union of the $\Gal(E/\Q)$-orbits in $S$ that contain
an element $\chi$ with $q(\chi)\not=0$. Since $\Gal(E/\Q)$ contains complex conjugation,
this condition is equivalent to the existence of $\chi$ such that $q(\chi)=1$.

The transcendental lattice $T(F_\C)$
is the smallest saturated sublattice of $P$ such that $P^{1,-1}\subset T(F_\C)\otimes\C$.
Thus $T(F_\C)=P\cap (T(F_\C)\otimes\Q)$, where $T(F_\C)\otimes\Q$ is 
the smallest $\Q$-vector space $V$
such that $\oplus_{q(\chi)=1}V_\chi\subset V\otimes E$. It is clear that $T(F_\C)\otimes E$ is the direct
sum of $V_\chi$ for $\chi\in S_\sharp$.

\ble \label{21jan2} The action of $\tau(1)$ on the $E$-vector space $P\otimes E$ is anti-linear
and preserves the Hodge decomposition.
We have $\tau(1)(\alpha_\chi)=\alpha_\chi$ for every $\chi\in\widehat G$. 
\ele
{\em Proof.} We have seen that the action of $\tau(1)$ on $\H^2(F_\C,\C(1))$ is anti-linear,
hence the action on $P\otimes E\subset \H^2(F_\C,\C(1))$ is also anti-linear.
The generator $2\pi i e'$ of $P$ is fixed by $\tau(1)$, see (\ref{conj}).
We have $\tau(1) [g]=[g^{-1}]\tau(1)$ for $g\in G$ and $\tau(1)\e=\e^{-1}\tau(1)$, thus
$\tau(\alpha_\chi)=\alpha_\chi$. $\Box$

\medskip

That $\tau(1)$ is anti-linear and preserves the Hodge decomposition is a general fact.
In Deligne's notation $\tau=F_\infty$, but $F_\infty(\H^{p,q})=\H^{q,p}$ by \cite[0.2.5]{Del}. 
Our $\tau(1)$ is the composition of $F_\infty$
with the complex conjugation of the coefficients $\C$, hence $\tau(1)(\H^{p,q})=\H^{p,q}$, 
see \cite[Prop.~1.4, Cor.~1.6]{Del}.

\subsection{Recovering full cohomology from primitive cohomology}

As was noted in Section \ref{prim_coh}, we have 
$H=\Z\Lambda\oplus P$, where $\Lambda\in H$ is the class of a line in $F_\C$.
Let $c=L-d\Lambda\in P$, where $L$ is the hyperplane section class. 
Then $H$ is the subgroup of $\Q L\oplus (P\otimes\Q)$ generated by $\Z L\oplus P$
and $\frac{1}{d}(L-c)$. We have $\langle L,L\rangle=d$, and 
Proposition \ref{o1} explicitly describes
the restriction of the cup-product form to $P$.
Thus, in order to recover the full cohomology group $H=\H^2(F_\C,\Z(1))$ as a 
$G$-module together with its unimodular 
$G$-invariant cup-product bilinear form from the $G$-module $P$ and the restriction
of the cup-product form to $P$,
it is enough to find a formula for $c$. This can be done using 
calculations of intersection indices from the paper of A.~Degtyarev and I.~Shimada \cite{DS}. 
Recall that $\phi(x)$ and $\rho(x,y)$ were introduced in (\ref{phi&rho}).
It is immediate to check that in the ring $\Z[x,y]/(x^d-1,y^d-1)$ we have the identity \cite[Lemma 4.5]{DS}
$$(1-y)\phi(xy)=(1-x)(1-y)\rho(x,y).$$
Since $1-x$ is a unit in $\Q[x]/(\phi(x))$, in the ring $R_x\otimes R_y$ we have $\phi(xy)=(1-x)\rho(x,y)$.

\bpr \label{mai}
Write $u_0=(u_1 u_2 u_3)^{-1}$ so that $u_0u_1u_2u_3=1$. Define
$$c=(1-u_0)^{-1}\phi(u_0u_1)\cdot(1-u_2)^{-1}\phi(u_2u_3)=\rho(u_0,u_1)\rho(u_2,u_3)\in P.$$ 
Then $\frac{1}{d}(L-c)$ is the class of the line in $F_\C$ given by 
$x_1=\zeta_{2d} x_0,\, x_3=\zeta_{2d} x_2$.
The group $H\subset \Q L\oplus (P\otimes\Q)$ is generated by $\Z L\oplus P$ and $\frac{1}{d}(L-c)$.
\epr
{\em Proof.} Let $\Lambda\in H$ be the class of this line.
The lattice $H$ is non-degenerate, so it suffices to show that the images of $\frac{1}{d}(L-c)$ 
and $\Lambda$ in $H^*\subset (\Z L)^*\oplus P^*$ coincide.
The images of these elements in $(\Z L)^*$ are equal, so it remains to show that they are equal in $P^*$.

The $G$-module isomorphism $P\cong \Z[G]/I$ gives rise to the map which evaluates
the linear function given by the cup-product with $x\in H$ on the elements of $G$:
$$ \mathrm{ev}\colon H\lra P^*\lra  \Z[G],\quad\quad \mathrm{ev}(x)= \sum_{g\in G}\langle x,g\rangle g$$
 By \cite[p.~989]{DS}, we have $\mathrm{ev}(\Lambda)=\psi$, where
$\psi=(1-u_1)(1-u_3)\phi(u_2u_3)$.
It is clear that $\mathrm{ev}(L)=0$, so
it remains to show that $\mathrm{ev}(c)=-d\psi$.
For $h\in G$, we find using the $G$-invariance of the cup-product pairing that
$$\mathrm{ev}(h)=\sum_{g\in G}\langle h, g\rangle g
=\sum_{g\in G}\langle 1, h^{-1}g\rangle g=
\sum_{g\in G}\langle 1, g\rangle gh=\mathrm{ev}(1)h.$$
We have 
$\mathrm{ev}(1)=(e|e)=-(1-u_0)(1-u_1)(1-u_2)(1-u_3)$, by (\ref{ee}).
Hence, we get
\begin{eqnarray*}
-\mathrm{ev}(c)&=&(1-u_0)(1-u_1)(1-u_2)(1-u_3)\rho(u_0,u_1)\rho(u_2,u_3)\\
&=&(1-u_1)(1-u_3)\phi(u_0u_1)\phi(u_2u_3)\\
&=&(1-u_1)(1-u_3)\phi(u_2u_3)^2\\
&=&(1-u_1)(1-u_3)d\phi(u_2u_3)=d\psi.\quad \Box
 \end{eqnarray*}

\subsection{Galois representation over the cyclotomic field}

Let $\ell$ be a prime number. Using the comparison theorem of classical and $\ell$-adic
\'etale cohomology, we obtain a $G$-module isomorphism
$$P_\ell\coloneqq P\otimes\Z_\ell=P^2(\ov F,\Z_\ell(1))\cong\Z_\ell[G]/(I\otimes\Z_\ell),$$ 
where $P^2(\ov F,\Z_\ell(1))\subset \H^2_\et(\ov F,\Z_\ell(1))$ 
is the orthogonal complement to the hyperplane section class.
Since $T_\ell\coloneqq T(F_\C)\otimes\Z_\ell$ is the orthogonal complement to the subgroup generated by
the classes of divisors in $P_\ell$, the Galois group $\Ga$ acts on $T_\ell$.

Let us assume that the ground field $k$ contains $E=\Q(\mu_d)$. 
Then the action of $\Ga=\Gal(\ov k/k)$ on $P^2(\ov F,\Z_\ell(1))$ commutes with the action of $G$.

Let $\lambda$ be a prime of $E$ above $\ell$. Let $O_\lambda$ be the ring of integers of $E_\lambda$.
The $E_\lambda$-vector space $P\otimes E_\lambda$ is the direct sum of the 
1-dimensional $G$-eigenspaces $V_\chi\otimes_E E_\lambda$, where $\chi\in S$, which are therefore 
preserved by the Galois group $\Ga$. We have 
$$T(F_\C)\otimes E_\lambda=\bigoplus_{\chi\in S_\sharp}V_\chi\otimes_E E_\lambda.$$

Let $O_k$ be the ring of integers of $k$.
Let $\mathfrak p\subset O_k$ be a prime ideal coprime to $d\ell$, and 
let $\F_{\mathfrak p}=O_k/\mathfrak p$. Let $p={\rm char}(\F_{\mathfrak p})$.
For $x\in\F_{\mathfrak p}^\times$ let $\psi(x)$ be the unique $d$-th root of unity in $E\subset k$ whose reduction
in $\F_{\mathfrak p}$ equals $x^{(N({\mathfrak p})-1)/d}$, 
where $N({\mathfrak p})=|\F_{\mathfrak p}|$. Thus $\psi\colon\F_{\mathfrak p}^\times\to\mu_d$
is a multiplicative character of order $d$. The Gauss sum $g(r)\in\Q(\e,\zeta)$ is 
$$g(r)=\sum_{x\in \F_{\mathfrak p}^\times}\psi(x)^r\,\zeta^{\Tr_{\F_{\mathfrak p}/\F_p}(x)}, 
\quad\quad r=0,\ldots,d-1,$$
where $\zeta$ is a primitive $p$-th root of unity in $\C$.
We have $g(0)=0$, $\ov{g(r)}=\psi(-1) g(-r)$ and $|g(r)|^2=N({\mathfrak p})$. Note that
$\psi(-1)=1$ if $d$ is odd or both $d$ and $(N({\mathfrak p})-1)/d$ are even;
if $d$ is even and $(N({\mathfrak p})-1)/d$ is odd, then $\psi(-1)=-1$, see \cite[\S 5.2]{LN}
or \cite[\S 10.3]{IR}.

The Jacobi sum $J(\chi)$ is an element of the ring of integers $O_E=\Z[\e]$ defined by
$$J(\chi)=\sum_{x_1+x_2+x_3=1}\psi(x_1)^l\psi(x_2)^m\psi(x_3)^n=\frac{g(l)g(m)g(n)}{g(l+m+n)}
=\psi(-1)\frac{g(l)g(m)g(n)g(r)}{N({\mathfrak p})},$$
where $x_1, x_2, x_3\in\F_{\mathfrak p}$.
The integer $r$ satisfies $l+m+n+r=\left \lceil{\frac{l+m+n}{d}}\right \rceil d$.
The second equality here is \cite[Thm. 5.21]{LN}, and the third one follows from the properties of Gauss sums
discussed above, cf. \cite[Ch.~8, Thm.~3]{IR}.

The following proposition can be found in D. Ulmer's paper \cite{U}, though it goes back to N. Katz and T. Shioda. 
We give a simpler proof using the Fourier transform on the finite abelian group $G$.

\bpr \label{o3}
Let $\lambda$ be a prime of $E=\Q(\mu_d)$ above $\ell$. 
Let $\mathfrak p$ be a prime of $k$ not dividing $d\ell$.
Then $\Frob_{\mathfrak p}$ acts on $V_\chi\otimes_E E_\lambda$, where $\chi\in S$, as multiplication by 
$$\psi(-1)N({\mathfrak p})^{-1}J(\chi)\ \in \ N({\mathfrak p})^{-1}O_E\subset E.$$
\epr

Before giving a proof we explain how the Galois representation attached to
an arbitrary diagonal surface can be obtained from the Galois representation attached to
the Fermat surface.

\brem \label{twisting general}
{\rm Let $K$ be a subfield of $\C$. For $a=(a_1,a_2,a_3)\in(K^\times)^3$
we write $X=X_a$ for the surface in $\P^3_K$ given by
$$x_0^d+a_1x_1^d+a_2x_2^d+a_3x_3^d=0.$$
In particular, the Fermat surface $F=X_{(1,1,1)}$.
The surface $X_a$ is obtained by twisting $F$ over $K$
by a 1-cocycle of $\Gal(\ov K/K)$ with values 
in the group $K$-scheme $G=(\mu_d)^3$, 
whose class in $\H^1(K,G)\cong(K^\times/K^{\times d})^3$ is the image of $a=(a_1,a_2,a_3)$
under the natural coordinate-wise map $(K^\times)^3\to (K^\times/K^{\times d})^3$.
The $\Gal(\ov K/K)$-modules $\Pic(\ov X_a)$ and $T_\ell(X_{a,\C})$ are obtained from 
the $\Gal(\ov K/K)$-modules $\Pic(\ov F)$ and $T_\ell(F_\C)$, respectively, by twisting 
them by this cocycle with respect to the relevant induced actions of $G$.
If $E\subset K$, then the action of $\Gal(\ov K/K)$ commutes with 
the action of $G$, so $\Gal(\ov K/K)$ preserves the eigenspaces 
$V_\chi\otimes_E E_\lambda$. In the case of $X_a$, the action of
$\Gal(\ov K/K)$ on $V_\chi\otimes_E E_\lambda$
is the same as the action in the case of $F$ followed by multiplication by 
the inverse of the $d$-power character $\Gal(\ov K/K)\to\mu_d$ associated to $a_1^la_2^ma_3^n$, where $\chi=(l,m,n)$.}
\erem

\noindent{\em Proof of Proposition \ref{o3}.} 
For $\chi=(l,m,n)\in S$ let
$h(\chi)\in E_\lambda$ be the eigenvalue of $\Frob_{\mathfrak p}$ on 
$V_\chi\otimes E_\lambda$.
If $\chi\in\widehat G$ is not in $S$, then we set $h(\chi)=0$.

Let $a=(a_1,a_2,a_3)\in G$. The multiplicative character $\psi$
defines an isomorphism $\F_{\mathfrak p}^\times/\F_{\mathfrak p}^{\times d}\tilde\lra\mu_d$. 
Pick any lifting $\tilde a_i\in \F_{\mathfrak p}^\times$ of $a_i$, for $i=1,2,3$, and consider
the surface in $X_a\subset\P^3_{\F_{\mathfrak p}}$ with equation
$$ x^d+\tilde a_1x_1^d+\tilde a_2x_2^d+\tilde a_3x_3^d=0.$$
Since $G$ acts on $V_\chi$ via $\chi$, the twisted action of $\Frob_{\mathfrak p}$ on $V_\chi\otimes E_\lambda$,
i.e.~the action defined by the action of the Galois group $\Ga$ on $P^2(\ov X_a,E_\lambda(1))$,
is via multiplication by $h(\chi)/\psi(\tilde a_1)^{l} \psi(\tilde a_2)^{m}\psi(\tilde a_3)^{n} $.
Hence, by the Lefschetz trace formula, we have 
$|X_a(\F_{\mathfrak p})|=|\P^2(\F_{\mathfrak p})|+f(a)N({\mathfrak p})$, where
$$f(a)=\sum_{\chi\in\widehat G} \psi(\tilde a_1^l \tilde a_2^m \tilde a_3^n)^{-1} h(\chi)
=\sum_{\chi\in\widehat G} \chi(-a) h(\chi). $$
Weil's classical calculation with exponential sums \cite{W1}, \cite{LN} gives that this
also holds when $h(\chi)$ is replaced by $g(\chi)=\psi(-1)J(\chi)N({\mathfrak p})^{-1}$
for $\chi\in S$ and by $g(\chi)=0$ for $\chi\notin S$.
Thus $f$ is the Fourier transform of either $h$ or $g$. The inverse Fourier transform now gives
$$h(\chi)=d^{-3}\sum_{a\in G}\chi(a)f(a)=g(\chi).\eqno{\Box}$$

\section{The Fermat quartic surface}

From now on let $d=4$.
In this section we let $k=E=\Q(i)$ and $\O=\Z[i]$. 
For a number field $K\subset\ov\Q$ we write $\Ga_K=\Gal(\ov\Q/K)$.

\subsection{Explicit transcendental lattice}

By definition, the transcendental lattice $T(F_\C)$ is 
the smallest saturated $\Z$-submodule of $P$
whose complexification contains $P^{1,-1}$. Thus we have
$$T(F_\C)=P\cap (V_{(1,1,1)}\oplus V_{(3,3,3)})\subset P\otimes k.$$

\ble \label{3.1}
We have $T(F_\C)=\Z w_1\oplus\Z w_2$, where
$$w_1=8(\alpha_{(1,1,1)}+\alpha_{(3,3,3)}), \quad 
w_2=8i(\alpha_{(1,1,1)}-\alpha_{(3,3,3)}).$$
We have 
\begin{equation}
\langle w_1,w_1\rangle=\langle w_2,w_2\rangle=8, \quad \langle w_1,w_2\rangle=0.
\label{ww}
\end{equation}
\ele
{\em Proof.} It is clear that $w_1, w_2$ is a basis of $T(F_\C)\otimes\Q$.
One immediately obtains (\ref{ww}) from (\ref{cup}). Recall
that $2^6\alpha_{(1,1,1)}$ is the sum of $i^{-(a+b+c)}u_1^au_2^bu_3^c$,
where $a,b,c\in\{0,1,2,3\}$. We calculate
$$\langle 2^6\alpha_{(1,1,1)}, u_1^\alpha u_2^\beta u_3^\gamma\rangle
=\sum_{a,b,c=0}^3 i^{-(a+b+c)}\langle u_1^{a-\alpha}u_2^{b-\beta}u_3^{c-\gamma},1\rangle=
4i^{-(\alpha+\beta+\gamma)},$$ 
hence
$\langle w_1,u_1^\alpha u_2^\beta u_3^\gamma\rangle={\rm Re}(i^{-(\alpha+\beta+\gamma)})$ and
$\langle w_2,u_1^\alpha u_2^\beta u_3^\gamma\rangle=-{\rm Im}(i^{-(\alpha+\beta+\gamma)})$.
This implies that if $sw_1+tw_2\in H$ for some $s,t\in \Q$, then $s,t\in\Z$.
This also implies $\langle w_i,P\rangle=\Z$. Recall that
$H=P\oplus\Z \Lambda$, where $\Lambda$ is the class of a line in $F_\C$.
Since $w_i\in T(F_\C)\otimes\Q$, we have
$\langle w_i,\Lambda\rangle=0$ for $i=1,2$, hence $\langle w_i, H\rangle=\Z$. 
This implies $w_i\in H^*=H$ by the unimodularity of $H$.
This finishes the proof that $T(F_\C)=\Z w_1\oplus\Z w_2$. $\Box$

\medskip

This lemma implies that $\Delta\simeq(\Z/8)^2$, as was already discovered in \cite{PSS}.

Let $\mu_4$ be one of the factors of $G=(\mu_4)^3$. In other words, choose one
of $x_1, x_2, x_3$, say $x_3$, and multiply $x_3$ by fourth roots of unity.
Let us write the induced action of $u\in\mu_4$ on $P$ as $[u]$. Extending
the coefficients, this gives rise to a $k$-linear action of $\mu_4$ on $P\otimes k$. By (\ref{29dec1})
we then have $[i]\alpha_{(l,m,n)}=i^n \alpha_{(l,m,n)}$. In particular,
$[i]\alpha_{(1,1,1)}=i \alpha_{(1,1,1)}$ and $[i]\alpha_{(3,3,3)}=-i \alpha_{(3,3,3)}$.
(Note that these formulae will be the same if we consider the action of $\mu_4$
on $x_1$ or $x_2$ instead of $x_3$.)
Hence we obtain that $[i]$ acts on $T(F_\C)$ by the rule 
$[i]w_1=w_2$, $[i]w_2=-w_1$.
This gives an isomorphism of $\O$-modules, where $\O=\Z[i]$:
\begin{equation}
T(F_\C)\cong\O. \label{O}
\end{equation}
The isomorphism is unique up to the multiplication by an element of $\mu_4$;
we fix the isomorphism by identifying $1\in\O$ with $w_1$.
The cup-product form on $T(F_\C)$ then becomes $\langle x,y\rangle=4\Tr_{k/\Q}(x\ov y)$.

The normalised cup-product form $\frac{1}{8}\langle x,y\rangle$ induces an isomorphism
$T(F_\C)\tilde\lra T(F_\C)^*$. By the compatibility between classical and $\ell$-adic \'etale 
cohomology, for any prime $\ell$ the normalised cup-product in $\ell$-adic \'etale 
cohomology $\frac{1}{8}\langle x,y\rangle_\ell$ gives an isomorphism of $\Ga$-modules 
$$T(\ov X)_\ell\cong T(F_\C)\otimes\Z_\ell\tilde\lra T(F_\C)^*\otimes\Z_\ell,$$ 
cf. Section \ref{cohtools}.
Thus from (\ref{e1}) and (\ref{e19}) we obtain isomorphisms of $\Ga$-modules
\begin{equation}
T(X_\C)\otimes\Q/\Z \cong \Br(\ov X) \quad\text{and}\quad
T(X_\C)/\ell^m \cong \Br(\ov X)[\ell^m], \label{e99}
\end{equation}
for any prime $\ell$ and any positive integer $m$. Since $T(F_\C)^*=\frac{1}{8}T(F_\C)$, the 
exact sequence
$(\ref{three})$ becomes the exact sequence of $\Ga$-modules
\begin{equation}
0\lra \Br(\ov F)[8]\lra\Br(\ov F)\stackrel{[8]}\lra\Br(\ov F)\lra 0. \label{x2}
\end{equation}

\subsection{Explicit Galois representation} \label{S3.1}

Let $\ell$ be a prime number and let $\lambda$ be a prime of $k=\Q(i)$ above $\ell$. 
Let $\mathfrak p$ be a prime of $k$ not dividing $2\ell$. The principal ideal
$\mathfrak p\subset\O$ has a unique generator $\pi$ which is a {\em primary prime} 
$\pi\equiv 1\bmod (1+i)^3$. (Note that this works both in the split and the inert cases.)

Pinch and Swinnerton-Dyer \cite{PSD}, using Weil \cite{W1,W2} who relied on Stickelberger, 
computed the eigenvalues of $\Frob_{\mathfrak p}$ 
acting on $P^2(\ov F,k_\lambda(1))=P\otimes k_\lambda$, where $F$
is the quartic Fermat surface.
The Galois group $\Ga_k$ commutes with $G$,
hence the representation of $\Ga_k$ in $P\otimes k_\lambda$ attached to $F$ 
is a direct sum
of 1-dimensional representations $V_\chi\otimes k_\lambda$ described in Proposition \ref{o3}. 
For $d=4$ the multiplicative character $\psi(x)$ is the biquadratic residue 
$\left(\frac{x}{\pi}\right)_4$.

Below, for symmetry reasons, we represent characters
$\chi\in\widehat G$ by 4-tuples of integers $(l,m,n,r)\in\{1,2,3\}^4$ with sum divisible by 4.

\begin{itemize}
\item $\chi=(1,1,1,1)$. By Lemma \ref{o2} we have $P^{1,-1}=V_{(1,1,1,1)}\otimes_k\C$.
By Proposition~\ref{o3} and the quartic Jacobi sum calculation \cite[Prop.~9.9.5]{IR},
the representation of $\Ga_k$ in $V_{(1,1,1,1)}\otimes_k k_\lambda$
sends $\Frob_{\mathfrak p}$ to $\pi/\ov\pi$. (In the inert case this equals~$1$.)
\item $\chi=(3,3,3,3)$. We have $P^{-1,1}=\ov{P^{1,-1}}=V_{(3,3,3,3)}\otimes_k\C$,
so the representation of $\Ga_k$ in $V_{(3,3,3,3)}\otimes_k k_\lambda$ is conjugate to the
representation of $\Ga_k$ in $V_{(1,1,1,1)}\otimes_k k_\lambda$.
\item The spaces $V_\chi$ for all other characters belong to $P^{0,0}$.
The space $V_{(2,2,2,2)}$ is defined over $\Q$. There are six spaces like
$V_{(1,1,3,3)}$ (obtained by permuting the coordinates of $\chi$) that come in conjugate pairs.
The action of $\Ga_k$ on all of these spaces is trivial. Finally, there are twelve spaces like
$V_{(1,2,2,3)}$; the action of $\Frob_{\mathfrak p}$ on them is via multiplication by $\psi(-1)$.
In particular, we recover the known fact that $\Ga_{\Q(\mu_8)}$ acts trivially on $\Pic(\ov F)$.
\end{itemize}

It is more delicate to make explicit the Galois representation over $\Q$.
In the case $d=4$ we use the fact that
$\Ga_\Q$ is the semi-direct product of $\Ga_k$ and 
$\Gal(\C/\R)\cong\Gal(\ov\Q/(\ov\Q\cap\R))\cong\Z/2$. 

For odd primes $\ell$
the following statement was obtained in \cite{IS} by a different method
(see the proof of \cite[Lemma 4.2]{IS}).

\bpr \label{Gal}
The representation of $\Ga_\Q=\Ga_k\rtimes \Gal(\ov \Q/\ov \Q\cap\R)$ in 
$T(F_\C)\otimes\Z_\ell\cong\O\otimes\Z_\ell$ is as follows:

if $\mathfrak p$ is coprime to $2\ell$, then $\Frob_{\mathfrak p}\in \Ga_k$ acts as
multiplication by $\pi/\ov \pi$ where $\pi$ is a primary generator of $\mathfrak{p}$; 

the generator of $\Gal(\ov \Q/\ov \Q\cap\R)$ 
acts on $T(F_\C)$ as $\tau(1)$, namely, $\tau(1)(w_1)=w_1$, $\tau(1)(w_2)=-w_2$, i.e.,
as the usual conjugation of $\O$.
\epr
{\em Proof.} 
Recall that $T(F_\C)=\Z w_1\oplus \Z w_2$ was made into an $\O$-module using the action of $[i]$
which sends $w_1$ to $w_2$ and $w_2$ to $-w_1$.
We know that on $(\O\otimes\Z_\ell)\alpha_{(1,1,1)}$ the Frobenius element
$\Frob_{\mathfrak p}$, where $\mathfrak p$ is coprime to $2\ell$, 
acts as multiplication by $\pi/\ov\pi\in (\O\otimes\Z_\ell)^*$. One then
immediately checks that $\Frob_{\mathfrak p}$ acts on 
$T(X_\C)\otimes\Z_\ell=\Z_\ell w_1\oplus \Z_\ell w_2$ as the 
$(2\times 2)$-matrix given by multiplication by $\pi/\ov\pi$ on $\O\otimes\Z_\ell$.

Let us prove the second statement.
We have comparison isomorphisms
$$\H^2(F(\C),\Z/\ell^n)\tilde\lra \H^2_\et(F_\C,\Z/\ell^n)\tilde\lra \H^2_\et(\ov F,\Z/\ell^n)$$
which are functorial with respect to the action of $\Gal(\C/\R)=\Gal(\ov \Q/(\ov \Q\cap\R))$. 
Thus the action of $\Gal(\C/\R)$ on $\H^2(F(\C),\Z)$
induced by the complex conjugation of $F(\C)$, that is, the action sending
the generator of $\Gal(\C/\R)$ to $\tau$, is compatible with the representation
of $\Gal(\ov \Q/(\ov \Q\cap\R))$ in $\H^2_\et(\ov F,\Z_\ell)$, cf. \cite[0.3]{Del}.
It follows that the representation
of $\Gal(\ov \Q/(\ov \Q\cap\R))$ in $\H^2_\et(\ov F,\Z_\ell(1))$
sends the generator to $\tau(1)$.
The representation of $\Gal(\ov \Q/(\ov \Q\cap\R))$ in 
$T(F_\C)\otimes\Z_\ell\subset\H^2_\et(\ov F,\Z_\ell(1))$
extends to 
$T(F_\C)\otimes k_\lambda\subset\H^2_\et(\ov F,k_\lambda(1))$.
The $k$-vector space $T(F_\C)\otimes k$ has basis $\alpha_{(1,1,1)}$ and 
$\alpha_{(3,3,3)}$, which are fixed by $\tau(1)$, see Lemma
\ref{21jan2}. From this we obtain the action of $\tau(1)$ on $w_1$ and $w_2$. $\Box$

\bpr \label{c0}
{\rm (i)} The $\Ga_k$-invariant subgroup of $\Br(\ov F)\{2\}$ is $\Br(\ov F)[4]$.

{\rm (ii)} The $\Ga_{k(\sqrt{2})}$-invariant subgroup
of $\Br(\ov F)\{2\}$ is $\Br(\ov F)[8]$.

{\rm (iii)} We have $k(\sqrt[4]{2})=k(\sqrt[4]{-2})$.
The $\Ga_{k(\sqrt[4]{-2})}$-invariant subgroup of $\Br(\ov F)\{2\}$ is
$\Br(\ov F)[16]$.

{\rm (iv)} The $\Ga_{k(\sqrt[8]{-2})}$-invariant subgroup of $\Br(\ov F)\{2\}$ is
$\Br(\ov F)[32]$.
\epr
{\em Proof.} (i) By (\ref{O}) and (\ref{e99}) there is
an isomorphism of $\Ga_k$-modules $\Br(\ov F)[2^m]=\O/2^m$ for each $m\geq 1$.
We have $\pi=a+bi$, where $a$ is odd and $b$ is even. Hence $\pi\equiv\ov\pi\bmod 4$. 
From the explicit description of the Galois action in Proposition \ref{Gal} we see that $\Ga_k$
acts trivially on $\Br(\ov F)[4]$. To finish the proof of (i) it is enough to show that
each element of $\O/8$ fixed by $\Ga_k$ is divisible by 2. One easily
checks that if $x+iy\in\O$ is fixed modulo 8 by $\Frob_\pi$, where $\pi=-1+2i$,
then $x$ and $y$ are even.

(ii) A supplement to the quadratic reciprocity law states that $2$ is a square 
modulo a prime $p\equiv 1\bmod 4$ precisely when $p\equiv 1\bmod 8$. 
Write $p=N(\pi)$, where $\pi$ is a primary prime in $\O$. 
Since $\Z/p\tilde\lra \O/\pi$, we see that $2$ is a square modulo $\pi$ if and
only if $p\equiv 1\bmod 8$, that is, if and only if $4$ divides $b$, that is,
if and only if $\pi\equiv\ov\pi\bmod 8$. This shows that $\Ga_{k(\sqrt{2})}$
acts trivially on $\Br(\ov F)[8]$. The same calculation as in (i) shows that
if $x+iy\in\O$ is fixed modulo 16 by $\Frob_\pi$, where $\pi=1+4i$,
then $x$ and $y$ are even.

(iii) A supplement to the biquadratic reciprocity law \cite[Thm.~4.23]{Cox},
\cite[Thm.~6.9]{lemmermeyer} states
that $2$ is the fourth power modulo $\pi$ if and only if $8$ divides $b$, that is,
if and only if $\pi\equiv\ov\pi\bmod 16$. The same calculation as in $(i)$ shows that
if $x+iy\in\O$ is fixed modulo $32$ by $\Frob_\pi$, where $\pi=5+8i$,
then $x$ and $y$ are even.

(iv) A supplement to the octic reciprocity law \cite[II.12]{western} states that for 
$\ov\pi\pi=p\equiv 1\bmod 8$ with $8\mid b$, $2$ is an eighth power modulo $\pi$ 
if and only if $b\equiv p-1\equiv 0\bmod 16$ or $b\equiv p-1\equiv 8\bmod 16$ 
and $-1$ is an eighth power modulo $\pi$ if and only if $p-1\equiv 0\bmod 16$. 
Thus $-2$ is an eighth power modulo $\pi$ if and only if $b\equiv 0\bmod 16$, that is $\pi\equiv\ov\pi\bmod 32$. 
The same calculation as in $(i)$ shows that
if $x+iy\in\O$ is fixed modulo $64$ by $\Frob_\pi$, where $\pi=1+16i$,
then $x$ and $y$ are even. $\Box$

\brem{\rm The invariant subfield of the kernel of the Galois action on
$\Br(\ov F)[2^n]$ is the ring class field of $k=\Q(i)$ of conductor $2^{n-1}$.
This follows from the definition: $K_c/K$ is the abelian extension of $K$ 
corresponding to the class group $I_K(c)/P_{K,\Z}$; in other words, precisely the primes
congruent to an element of $\Z$ modulo $c$ split in $K_c$. Indeed,
$\Frob_\pi$ acts on $\O/2^n$ trivially if and only if $\pi\equiv\ov\pi\bmod 2^n$,
which is equivalent to the divisibility of $b$ by $2^{n-1}$. See \cite[p.~179]{Cox}.
For the calculation $\Q(i)_8=\Q(i,\sqrt[4]{2})$, see \cite[Prop.~9.5]{Cox}.}
\erem

We now reprove a result of Ieronymou \cite[Thm.~4.2]{I}. Recall that 
if $K_1\subset K_2$ are algebraically closed fields of characteristic zero,
and $X$ is a smooth and proper variety over $K_1$, 
then the induced map $\Br(X)\to \Br(X\times_{K_1}K_2)$ is an isomorphism.
Thus we can use the notation $\Br(\ov F)$ without specifying
the underlying algebraically close field.

Since $F$ is defined over $\Q$, 
from the exact sequence (\ref{x2}) we obtain
the exact sequence of $\Ga_\Q$-modules
\begin{equation}
0\lra \O/8\stackrel{[2]}\lra\O/16\lra\O/2=\Br(\ov F)[2]\lra 0.\label{23jan1}
\end{equation}

The following corollary should be compared to Ieronymou's results \cite[Thm. 3.1, Prop. 4.9]{I}.

\bco  \label{c11}
Let $F$ be the Fermat quartic surface.
Let $K$ be a field extension of $\Q(i,\sqrt[4]{2})$.
Then the natural map $\Br(F_K)[2]\to\Br(\ov F)[2]$ is surjective.
\eco
{\em Proof.} It is enough to prove this for $K=\Q(i,\sqrt[4]{2})$.
The Fermat quartic surface has obvious rational points over $\Q(\mu_8)$;
pick one of them and call it $P$. We can write $\Br(F_K)=\Br(K)\oplus \Br(F_K)_P$,
where $\Br(F_K)_P\subset \Br(F_K)$ consists of the classes vanishing at $P$. 
Since $\Ga_K=\Gal(\ov\Q/K)$ acts trivially on $\Pic(\ov F)\simeq\Z^{20}$,
we have $\H^1(K,\Pic(\ov F))=0$. The standard spectral sequence then shows that
the canonical map $\Br(F_K)\to\Br(\ov F)$ restricted to $\Br(F_K)_P$
is an injective map $\Br(F_K)_P\hookrightarrow \Br(\ov F)$.
It follows that an element of $\Br(\ov F)[2]$
which can be lifted to an element of $\Br(F_K)$, can also be lifted to an element of $\Br(F_K)[2]_P
\subset \Br(F_K)[2]$.

It remains to show that $\Br(\ov F)[2]$ is contained in the image of $\Br(F_K)$ in $\Br(\ov F)$.
By Proposition \ref{c0} (iii) the abelian groups in (\ref{23jan1}) are trivial 
$\Ga_K$-modules. 
It follows that $\partial_1(\Br(\ov F)[2])=0$, so we can conclude by applying
Corollary \ref{1}. $\Box$

\subsection{Diagonal quartic surfaces}

We start with a general observation about arbitrary diagonal quartic surfaces
(proved in \cite[Prop.~4.8]{Sko17} by a different method).

\bpr \label{nonzero}
Let $X=X_a$ be a diagonal quartic surface over a field $K$ of characteristic zero.
We have $\Delta_X[2]^{\Gal(\ov K/K)}\not=0$, hence
$\Br(\ov X)[2]^{\Gal(\ov K/K)}\not=0$.
\epr
{\em Proof.} Without loss of generality we assume that $K=\Q(a_1,a_2,a_3)$;
in particular, $K$ can be embedded into $\C$. Then we have a canonical isomorphism
$\Pic(\ov X)\cong\Pic(X_\C)$.

We construct an explicit $\Gal(\ov K/K)$-invariant 2-torsion element in $\Delta_X$
using the canonical isomorphism of $\Gal(\ov K/K)$-modules
$\Delta_X\cong\Pic(\ov X)^*/\Pic(\ov X)$.
The snake lemma shows that the $\Gal(\ov K/K)$-module
$\Delta_X[2]$ is canonically isomorphic to the kernel
of the map $\Pic(\ov X)/2\to \Pic(\ov X)^*/2$ given by the intersection pairing. Thus it is enough
to exhibit a divisor class $D\in\Pic(\ov X)$ with the following properties:

(a) $\langle D,x\rangle\in2\Z$ for every $x\in \Pic(\ov X)$;

(b) the class of $D$ in $\Pic(\ov X)/2$ is $\Gal(\ov K/K)$-invariant; 

(c) the class of $D$ in $\Pic(\ov X)/2$ is not zero.

\noindent Consider the surjective morphism $s:X\to Q$,
where $Q$ is the diagonal quadric $$y_0^2+a_1y_1^2+a_2y_2^2+a_3y_3^2=0,$$ and
$s\big((x_0:x_1:x_2:x_3)\big)=(x_0^2:x_1^2:x_2^2:x_3^2)$. There are 
two families of lines on $Q$ defined over $\ov K$.
Let $f_\pm\in\Pic(\ov Q)$ be the classes of these families, and let $D_\pm=s^*(f_\pm)$.
It is clear that the set $\{D_+,D_-\}\subset \Pic(\ov X)$ is $\Gal(\ov K/K)$-invariant.

Let us show that we can take $D$ to be either of $D_+$ and $D_-$. Indeed,
the inverse images of the two families of lines in $Q$ are pencils of genus 1 curves on $F$
such that $D_++D_-=2L$,
where $L$ is the hyperplane section class.
It is well known and easy to check
that each of the 48 lines on $\ov X$ is a component of a degenerate fibre of one of these families.
Thus for any line $\ell\subset \ov X$ we have either $\langle\ell, D_+\rangle=0$ or 
$\langle\ell, D_+\rangle=2$. Since $\Pic(\ov X)$ is generated by lines, we get 
$\langle D_+, x\rangle\in2\Z$ for every $x\in \Pic(\ov X)$. By the formula
$D_++D_-=2L$, the same holds for $D_-$. This gives (a) and (b).

The class of an integral curve of arithmetic genus 1 on a complex K3 surface $X_\C$ is
primitive in $\Pic(X_\C)=\Pic(\ov X)$ \cite[Ch.~2, Rem.~3.13 (i)]{Huybrechts}.
This gives (c). $\Box$

\bde
Let $X$ be a diagonal surface over a field $K$.
We say that $X$ is {\bf split} by a field extension $L/K$ if $X_L\cong F_L$.
\ede

\bthe \label{c1}
Let $X$ be a diagonal quartic surface over a field $K$ containing $k=\Q(i)$, which is
split by a field extension of $K$ not containing $\sqrt[4]{2}$.
Then $$\Br(\ov X)[2]\cap \Im[\Br(X)\to\Br(\ov X)]=0.$$
\ethe
{\em Proof.} To prove that this intersection is zero we can enlarge the field $K$. 
So we can assume that $K$ contains $\Q(\mu_8)=k(\sqrt{2})$ 
but does not contain $\sqrt[4]{2}$, and that $X\simeq F_K$. 
The group $\Gal(\ov K/K)$ acts
on the terms of the exact sequence (\ref{23jan1})
via the map $\Gal(\ov K/K)\to \Gal(\ov\Q/k(\sqrt{2}))=\Ga_{k(\sqrt{2})}$.
By Proposition \ref{c0} (iii), $\Ga_{k(\sqrt{2})}$ acts 
via the surjective map 
$\Ga_{k(\sqrt{2})}\to\Gal(k(\sqrt[4]{2})/k(\sqrt{2}))$. Since $K$ and $k(\sqrt[4]{2})$
are linearly disjoint over $k(\sqrt{2})$, the following map is surjective:
$$\Gal(\ov K/K)\to\Gal(K(\sqrt[4]{2})/K)\cong
\Gal(k(\sqrt[4]{2})/k(\sqrt{2}))\cong\Z/2.$$
Proposition \ref{c0} (ii) says that the $\Ga_{k(\sqrt{2})}$-invariant subgroup of $\O/16$ is $2\O/16$, which is exactly
the kernel of $\O/16\to\O/2$. Hence this is also
the $\Gal(\ov K/K)$-invariant subgroup of $\O/16$.
This implies that $\partial_1$ induces an injective map 
$$\partial_1\colon \Br(\ov F)[2]^{\Ga_K}=\O/2\lra \H^1(K,\Delta_F)=\H^1(K,\O/8).$$
We claim that in our situation 
$\partial_2 \colon \H^1(K,\Delta_F)\to\H^2(K,\Pic(\ov F))$ is injective.
Indeed, since $k(\sqrt{2})=\Q(\mu_8)$ is contained in $K$,
and $\Ga_{\Q(\mu_8)}$ acts trivially on $\Pic(\ov F)$,
the group $\Gal(\ov K/K)$ acts
trivially on $\Pic(\ov F)$, and hence also on $\Pic(\ov F)^*$.
Thus $\H^1(K,\Pic(\ov F)^*)=0$, so that $\partial_2$ is indeed injective.
We have proved that $\partial_2\partial_1$ is injective on $\Br(\ov F)[2]$,
so by Corollary \ref{1} the intersection of $\Br(\ov F)[2]$ 
and the image of $\Br(F_K)$ in $\Br(\ov F)$ is zero. $\Box$

\medskip

The following result shows that for any diagonal quartic surface over $\Q$
or $k=\Q(i)$ the Galois-invariant subgroup of $\Br(\ov X)$ is annihilated by $4$.

\bpr \label{c2}
Let $X=X_a$ be a diagonal quartic surface over $\Q$.

{\rm (i)} If $a_1a_2a_3\in k^{\times4}$, then the $\Ga_k$-invariant subgroup
of $\Br(\ov X)\{2\}$ is $\Br(\ov X)[4]$.

{\rm (ii)} If $a_1a_2a_3\in \Q^{\times4}$, then the $\Ga_{\Q}$-invariant subgroup
of $\Br(\ov X)\{2\}$ is $\frac{1}{4}(\Z+2i\Z)/\O$. If $a_1a_2a_3\in -4\Q^{\times4}$, 
then the $\Ga_{\Q}$-invariant subgroup
of $\Br(\ov X)\{2\}$ is $\frac{1}{4}(1-i)\Z/\Z$.

{\rm (iii)} If $a_1a_2a_3\in k^{\times2}$ but $a_1a_2a_3\notin k^{\times4}$, then 
the $\Ga_k$-invariant subgroup
of $\Br(\ov X)\{2\}$ is $\Br(\ov X)[2]$. The same holds if $k$ is replaced by $\Q$.

{\rm (iv)} If $a_1a_2a_3\notin k^{\times2}$, then the $\Ga_k$-invariant subgroup
of $\Br(\ov X)\{2\}$ is $\Z/2$ generated by the class of $1+i$ in $\frac{1}{2}\O/\O$.
\epr

\brem \label{twisting 4}
{\rm Let us recall the description of the twisted Galois representation from Remark \ref{twisting general}
in the case $d=4$. The group scheme $G=(\mu_4)^3$ acts on 
$\alpha_{(1,1,1)}$ via the product character $(\mu_4)^3\to\mu_4$.
If $K$ is a subfield of $\C$ and $a=(a_1,a_2,a_3)\in(K^\times)^3$,
then the action of $\Gal(\ov K/K)$ on $T(X_{a,\C})\otimes\Z_\ell$
is obtained by twisting the action of $\Gal(\ov K/K)$ on $T(F_\C)\otimes\Z_\ell$
by the 1-cocycle $\Gal(\ov K/K)\to\mu_4$ with class 
$a_1a_2a_3\in K^\times/K^{\times4}=\H^1(K,\mu_4)$.
For example, if $a_1a_2a_3\in K^{\times4}$, then 
the action of $\Gal(\ov K/K)$ on $T(X_{a,\C})\otimes\Z_\ell$
is the same as the action of $\Gal(\ov K/K)$ on $T(F_\C)\otimes\Z_\ell$.
In general, if $K$ contains $k$, the twisted action is obtained by 
composing the untwisted action with the inverse of 
the quartic character given by $a_1a_2a_3$.
The value of this character on $\Frob_\mathfrak p$ is 
$\left(\frac{a_1a_2a_3}{\pi}\right)_4^{-1}$.}
\erem

\noindent{\em Proof of Proposition \ref{c2}.} (i) By Remark \ref{twisting 4},
in this case the action of $\Ga_k$ is the same as for the
Fermat surface over $k$, so we read off the result from Proposition~\ref{c0}~(i).

(ii) By the same remark, the action $\Ga_\Q$ is the same as for the Fermat
surface, so compared to (a) we only need to take into account the action of $\tau(1)$. 
As explained in Proposition \ref{Gal},
it acts on $\O$ as the usual conjugation, hence the first formula.
If $a_1a_2a_3\in -4\Q^{\times4}$, then $\tau(1)$ acts on $\O$ as the usual complex conjugation followed by
multiplication by $(1+i)/(1-i)$, which the value on the complex conjugation of the 1-cocycle given by $-4$.
Thus the twisted action of $\tau(1)$ sends $a+bi$ to $b+ai$, hence the result.

(iii) Recall that $\Frob_{\mathfrak p}\in\Ga_k$ acts on
$\Br(\ov X)[2^m]=\O/2^m$, $m\geq 1$, as multiplication by 
$(\frac{a_1a_2a_3}{\mathfrak p})_4^{-1} \frac{\pi}{\bar\pi}$.
Here $(\frac{a_1a_2a_3}{\mathfrak p})_4$ is $1$ if $a_1a_2a_3$ is a fourth power modulo $\mathfrak p$,
and $-1$ otherwise. One easily checks that,
whatever the sign, the $\Frob_{-1+2i}$-invariant subgroup of $\O/2^m$ is killed by $4$.
But $\Frob_{\mathfrak p}$, for any $\mathfrak p$ coprime to $2$, acts on $\O/4$
as multiplication by $(\frac{a_1a_2a_3}{\mathfrak p})_4^{-1}$. Since $a_1a_2a_3$ is
not a fourth power in $k$, by the Chebotarev density theorem we can find a prime ideal $\mathfrak p$
such that $\Frob_{\mathfrak p}$ acts as $-1$. This proves (iii) for $k$.

The kernel of the natural map $\Q^\times/\Q^{\times4}\to k^\times/k^{\times4}$ is the group of order $2$
generated by $-4$. Hence a square in $\Q$ which is not a fourth power in $\Q$ cannot
become a fourth power in $k$. The action of $\tau(1)$ on $\O/2=\Br(\ov F)[2]$
is trivial, so the action of $\tau(1)$ on $\O/2=\Br(\ov X)[2]$ is by multiplication by the inverse of the value
of the 1-cocycle given by $a_1a_2a_3$ evaluated at the complex conjugation. 
This is an element of $\{\pm 1\}$, so it
acts trivially on $\O/2$. Thus the statement (iii) for $\Q$ follows from the same statement for $k$.

(iv) In this case $\Frob_{-1+2i}$ acts on $\O/2^m$ as $\frac{-1+2i}{-1-2i}=1+4i$ 
multiplied by $\pm 1$ or $\pm i$.
One checks that in all cases the $\Frob_{-1+2i}$-invariant subgroup of $\O/2^m$ is killed by $4$.
By the same argument as in (iii),
for any odd prime ideal $\mathfrak p\subset\O$, the action of $\Frob_{\mathfrak p}$ on $\O/4$ 
is by multiplication by $\pm 1$ or $\pm i$.
Since $a_1a_2a_3$ is not a square in $k$, by the Chebotarev density theorem we can find a
prime ideal $\mathfrak p$
such that $\Frob_\pi$ acts on $\O/4$ as multiplication by $\pm i$. Then the $\Frob_\pi$-invariant
subgroup of $\O/4$ is $\Z/2$ generated by $2(1+i)$. This subgroup is $\Ga_k$-invariant,
which proves (iv). $\Box$

\medskip

The following statement answers Question 2 raised at the end of \cite{Sko17} and gives a different 
proof of Proposition \ref{nonzero}.

\bco \label{coco}
Let $K$ be a field of characteristic zero, $a=(a_1,a_2,a_3)\in(K^\times)^3$, and let $X=X_a$.
Then the action of $\Gal(\ov K/K)$ on $\Br(\ov X)[2]$ fixes the unique non-zero $\mu_4$-invariant element;
the other two non-zero $2$-torsion elements are permuted by $\Gal(K(\sqrt{a_1a_2a_3})/K)$.
\eco
{\em Proof.} The Galois group $\Gal(\ov K/K)$ acts on $\Br(\ov F)$ via $\Ga=\Gal(\ov\Q/\Q)$,
but $\Ga$ acts trivially on $\Br(\ov F)[2]$ by the proof of (ii) above. Thus $\Gal(\ov K/K)$
acts on $\Br(\ov X)[2]$ via the inverse of the quartic character associated to $a_1a_2a_3$.
The action of $\mu_4$ on $\O/2\O$ fixes the class of $1+i$ and permutes the classes of $1$ and $i$,
hence the result. $\Box$

\medskip

Another result of Ieronymou \cite[Thm.~5.2]{I} is a corollary of Theorem \ref{c1}.

\bpr[Ieronymou] 
Let $X=X_a$ be a diagonal quartic surface over $\Q$ such that
$2$ does not belong to the subgroup of
$\Q^\times$ generated by $-4, a_1, a_2, a_3$ and $\Q^{\times4}$. Then 
$\Br(\ov X)[2]\cap \Im[\Br(X)\to\Br(\ov X)]=0$, so that
$\Br(X)[2]\subset\Br_1(X)$.
\epr
{\em Proof.} The surface $X$ is split by the field
$K=k(\sqrt[4]{a_1},\sqrt[4]{a_2},\sqrt[4]{a_3})$, where $k=\Q(i)$.
The kernel of 
$\Q^\times/\Q^{\times4}\to k^\times/k^{\times4}$
is generated by the class of $-4$, so the kernel of 
$\Q^\times/\Q^{\times4}\to K^\times/K^{\times4}$ is generated by the classes of $-4$, $a_1$, $a_2$, $a_3$.
Hence $2\notin K^{\times4}$, so the statement follows from Theorem \ref{c1}. $\Box$

\bpr \label{p1}
Let $X=X_a$ be a diagonal quartic surface over a field 
$K$ containing $k=\Q(i)$. 
If $a_1a_2a_3\in K^{\times4}$ and $2\in K^{\times4}$, then $\Br(\ov X)[2]\subset\Im[\Br(X)\to\Br(\ov X)]$.
\epr
{\em Proof.} By Remark \ref{twisting 4}, the action of $\Gal(\ov K/K)$ on $\Br(\ov X)$
is the same as the action of $\Gal(\ov K/K)$ on $\Br(\ov F)$, so the arguments from the 
second part of the proof of Corollary \ref{c11} still work. $\Box$

\medskip

The proof of Corollary \ref{c11} shows that in this situation
each element of $\Br(\ov X)[2]$ can be lifted to an element of $\Br(X)/\Br_0(X)$ of order $2$.

\section{Main theorem}

\subsection{Reduction to finite group cohomology}

Let $X_a$ be a diagonal quartic surface over $\Q$ with coefficients $a_1,a_2,a_3\in\Q^\times$.
To $a=(a_1,a_2,a_3)$ and a field $L$ of characteristic zero we attach a finite Galois extension 
$L_a$ and its Galois group
$$L_a=L(i,\sqrt[4]{2},\sqrt[4]{a_1},\sqrt[4]{a_2},\sqrt[4]{a_3}), \quad\quad
G_{a,L}=\Gal(L_a/L).$$
We shall also consider a larger Galois extension $\widetilde L_a$ and its Galois group
$$\widetilde L_a=L(i,\sqrt[8]{-2},\sqrt[4]{a_1},\sqrt[4]{a_2},\sqrt[4]{a_3}), \quad\quad
\widetilde G_{a,L}=\Gal(\widetilde L_a/L).$$
Note that $L_a$ splits $X_a$ and contains $\Q(\mu_8)$. We see that
$\Pic(\ov X)$, and hence also $\Pic(\ov X)^*$ and $\Delta_X$, are trivial $\Gal(\ov\Q/L_a)$-modules.
By Proposition \ref{c0} (iii), $\Br(\ov X)[16]$ is a trivial $\Gal(\ov\Q/L_a)$-module too.
Thus $\Ga_L=\Gal(\ov\Q/L)$ acts on these modules via $G_{a,L}$.
By Proposition \ref{c0} (iv), $\Br(\ov X)[32]$ is a trivial $\Gal(\ov\Q/\widetilde L_a)$-module,
so $\Ga_L$ acts on it via $\widetilde G_{a,L}$.

In view of Proposition \ref{c2}, the main theorem can be proved by
applying Corollary \ref{1.5} with $n=4$ to all possible groups $\widetilde G_{a,\Q}$ and $\widetilde G_{a,k}$,
where $k=\Q(i)$. This is a finite calculation because each of these groups is a subgroup of the 
`generic' Galois group $\widetilde G$ defined by treating $a_1$, $a_2$, $a_3$ as independent variables.
However, we have chosen a slightly more elaborate line of proof which has the advantage of reducing the use of
computer to a small number of group cohomology  calculations.

\medskip

The exact sequence (\ref{three}) in our particular case takes the form (\ref{x2}).
From it we obtain exact sequences of $\Ga_L$-modules
\begin{equation}
0\lra \Delta_X\lra\Br(\ov X)[16]\stackrel{[8]}\lra \Br(\ov X)[2]\lra 0, \label{eins}
\end{equation}
\begin{equation}
0\lra \Delta_X\lra\Br(\ov X)[32]\stackrel{[8]}\lra \Br(\ov X)[4]\lra 0. \label{zwei}
\end{equation}
The sequence (\ref{eins}) is thus an exact sequence of $G_{a,L}$-modules, and so
defines a boundary map $\H^0(G_{a,L},\Br(\ov X)[2])\to \H^1(G_{a,L},\Delta_X)$.
Similarly, the sequence (\ref{zwei}) defines a boundary map
$\H^0(\widetilde G_{a,L},\Br(\ov X)[4])\to \H^1(\widetilde G_{a,L},\Delta_X)$.

\ble \label{new}
Let $k=\Q(i)$.
Suppose that the image of the natural map $G_{a,L}\to\Gal(k(\sqrt[4]{2})/\Q)$
contains the subgroup $\Gal(k(\sqrt[4]{2})/k)\simeq\Z/4$. Then

{\rm (a)} the group $\Br(\ov X)[2]\cap{\rm Im}[\Br(X_L)\to\Br(\ov X)]$ is the intersection of images
of the following maps:
$$\H^1(G_{a,L},\Pic(\ov X)^*)\lra\H^1(G_{a,L},\Delta_X)\longleftarrow\H^0(G_{a,L},\Br(\ov X)[2]);$$

{\rm (b)} the group $\Br(\ov X)\{2\}\cap{\rm Im}[\Br(X_L)\to\Br(\ov X)]$ is the intersection of images
of the following maps:
$$\H^1(\widetilde G_{a,L},\Pic(\ov X)^*)\lra\H^1(\widetilde G_{a,L},\Delta_X)\longleftarrow
\H^0(\widetilde G_{a,L},\Br(\ov X)[4]).$$
\ele
{\em Proof.} Let us show that $\Br(X)\{2\}^{\Ga_L}$ is killed by 4. For this it is enough to show that
$\Br(X)[8]^{\Ga_L}$ is killed by 4. Let us first assume that $a_1a_2a_3$ is a fourth power
in $L$. In this case, the $\Ga_L$-modules $\Br(\ov X)\{2\}$ and $\Br(\ov F)\{2\}$ are isomorphic.
By Proposition \ref{c0} (ii), the group $\Ga_L$ acts on $\Br(\ov F)[8]$ via its image
$\Gal(L(i,\sqrt{2})/L)$. By assumption, the subgroup $\Gal(L(i,\sqrt{2})/L(i))$ is 
isomorphic to $\Gal(k(\sqrt{2})/k)\simeq\Z/2$. As in the proof of 
Proposition \ref{c0} (i), the non-trivial element of $\Gal(k(\sqrt{2})/k)$ can be taken to
be ${\rm Frob}_\pi$, where $\pi=-1+2i$. If $x+iy\in\O$ is fixed modulo 8
under the multiplication by $\frac{\pi}{\bar\pi}$, then $x$ and $y$ are even.
Thus 4 kills $\Br(\ov X)\{2\}^{\Ga_L}$ in this case. In general, the Galois representation
on $\Br(\ov X)[8]$ is obtained by twisting the representation on $\Br(\ov F)[8]$ 
by the inverse of the quartic character of $a_1a_2a_3$, so
$\frac{\pi}{\bar\pi}$ is multiplied by $\pm 1$ or $\pm i$. An easy check shows that 
4 kills $\Br(\ov X)\{2\}^{\Ga_L}$ in all cases.

In view of the exact sequence (\ref{eins}) this implies that the boundary maps 
$$\Br(\ov X)[2]^{\Ga_L}\lra \H^1(L,\Delta_X), \quad \quad
\H^0(G_{a,L},\Br(\ov X)[2])\lra \H^1(G_{a,L},\Delta_X)$$
are injective. We have a commutative diagram
$$\xymatrix{
\H^1(L,\Pic(\ov X)^*)\ar[r]&\H^1(L,\Delta_X)&
\Br(\ov X)[2]^{\Ga_L}\ar@{_{(}->}[l]_{\partial_1}\\
\H^1(G_{a,L},\Pic(\ov X)^*)\ar[r]\ar[u]_{\cong}&\H^1(G_{a,L},\Delta_X)\ar@{_{(}->}[u]&
\H^0(G_{a,L},\Br(\ov X)[2])\ar@{_{(}->}[l]\ar[u]_{\cong}
}$$
where the vertical arrows are inflation maps. Claim (a) follows from this diagram.

Similarly, the boundary maps defined by the exact sequence (\ref{zwei})
$$\Br(\ov X)\{2\}^{\Ga_L}=\Br(\ov X)[4]^{\Ga_L}\lra \H^1(L,\Delta_X), \quad
\H^0(\widetilde G_{a,L},\Br(\ov X)[4])\lra \H^1(G_{a,L},\Delta_X)$$
are injective. We have a similar commutative diagram
$$\xymatrix{
\H^1(L,\Pic(\ov X)^*)\ar[r]&\H^1(L,\Delta_X)&
\Br(\ov X)[4]^{\Ga_L}\ar@{_{(}->}[l]_{\partial_1}\\
\H^1(\widetilde G_{a,L},\Pic(\ov X)^*)\ar[r]\ar[u]_{\cong}&\H^1(\widetilde G_{a,L},\Delta_X)\ar@{_{(}->}[u]&
\H^0(\widetilde G_{a,L},\Br(\ov X)[4])\ar@{_{(}->}[l]\ar[u]_{\cong}.
}$$
It immediately implies (b). $\Box$

\subsection{The case of $\Q(i)$}

The main result of this section is the following description of the 2-primary torsion subgroup
of the transcendental Brauer group of $X$ over the field $k=\Q(i)$. 

\bthe \label{main1}
Let $X=X_a$ be a diagonal quartic surface with $a_1, a_2, a_3\in \Q^\times$.
Then $\left(\Br(X_k)/\Br_1(X_k)\right)\{2\}$ is $\Z/2$
if $X_k$ is equivalent to $X_{1,2,8}\times_\Q k$, and $0$ otherwise.
\ethe

We will apply Lemma \ref{new} to a field extension $L=K$ of $k$ 
not containing $k(\sqrt{2})$ and such that $K_a=K(\sqrt[4]{2})$.
Then
$$G_{a,K}=\Gal(K(\sqrt[4]{2})/K)=\Gal(k(\sqrt[4]{2})/k),$$
$$\widetilde G_{a,K}=\Gal(K(\sqrt[8]{-2})/K)=\Gal(k(\sqrt[8]{-2})/k).$$
In particular, the condition of Lemma \ref{new} is satisfied.

\ble \label{drei}
Let $K$ be a field extension of $k$ not containing $k(\sqrt{2})$ such that $K_a=K(\sqrt[4]{2})$.
Then the cardinality
of the image of $\H^1(G_{a,K},\Pic(\ov X)^*) \to \H^1(G_{a,K},\Delta_X)$ is equal to the cardinality
of the cokernel of $\H^0(G_{a,K},\Pic(\ov X)^*)\to \H^0(G_{a,K},\Delta_X)$. The same is true if $G_{a,K}$
is replaced by $\widetilde G_{a,K}$.
\ele
{\em Proof.} By the duality theorem for Tate cohomology groups
\cite[Ch.~VI, \S 7, Exercise 3]{Brown} we have an isomorphism
$$\H^{-1}(G_{a,K},\Pic(\ov X)^*)\cong\Hom(\H^1(G_{a,K},\Pic(\ov X)), \Z/4).$$
Since $G_{a,K}$ is cyclic, we have a canonical isomorphism 
$$\H^{-1}(G_{a,K},\Pic(\ov X)^*)\cong\H^1(G_{a,K},\Pic(\ov X)^*).$$
Thus $\H^1(G_{a,K},\Pic(\ov X)^*)$ is (non-canonically) isomorphic to $\H^1(G_{a,K},\Pic(\ov X))$.
The group $\widetilde G_{a,K}$ is also cyclic, so the same argument applies with $\Z/4$ replaced by $\Z/8$.

We have an exact sequence of abelian groups
$$
\H^0(G_{a,K},\Pic(\ov X)^*)\to \H^0(G_{a,K},\Delta_X)\to \H^1(G_{a,K},\Pic(\ov X))
\to \H^1(G_{a,K},\Pic(\ov X)^*), 
$$
and a similar sequence for $\widetilde G_{a,K}$.
The cardinalities of the third and fourth groups are equal, 
and our statement follows. $\Box$

\ble \label{good}
Let $K$ be a field extension of $k$ not containing $k(\sqrt{2})$.
Let $X=X_a$ be a diagonal quartic surface where $(a_1,a_2,a_3)$ is in the following list:
$$(1,1,2), \quad (1,1,4), \quad (1,1,8), \quad (1,2,2), \quad (1,2,4), \quad (1,4,4), \quad (1,4,8), \quad (2,4,8).$$
Then $\Br(\ov X)[2]\cap {\rm Im}[\Br(X_K)\to\Br(\ov X)]=0$.
\ele
{\em Proof.} In the proof of the statement, we can replace $K$ by a finite extension.
In the cases $(1,1,4)$ and $(1,4,4)$ we go over to $K(\sqrt{2})$ which splits $X$ but does not contain
$\sqrt[4]{2}$, so the result follows from Theorem \ref{c1}.

We identify $(1,1,2)$ and $(1,1,8)$ as case B59, while $(2,4,8)$
is case B46 in \cite[Table B]{Bright}. In these cases $\H^1(G_{a,K},\Pic(\ov X))=0$, 
hence $\H^1(G_{a,K},\Pic(\ov X)^*)=0$, so we conclude by Lemma
\ref{new} (a).

It remains to consider the cases $(1,2,2)$, $(1,2,4)$, and $(1,4,8)$. 
Using a computer we check that the map 
$\H^0(G_{a,K},\Pic(\ov X)^*)\to \H^0(G_{a,K},\Delta_X)$ is surjective
in these three cases (see Lemma~\ref{app-surj}). By Lemma \ref{drei} it follows that
$\H^1(G_{a,K},\Pic(\ov X)^*) \to \H^1(G_{a,K},\Delta_X)$ is the zero map.
We conclude by applying Lemma \ref{new} (a). $\Box$

\medskip

\noindent {\bf Proof of Theorem} \ref{main1}. 
Since $-4$ is a fourth power in $k=\Q(i)$ we can multiply the coefficients $a_1, a_2, a_3$ 
by $-4$ and fourth powers of their denominators
if necessary and so assume that they are positive integers. Write
$$a_1=2^\alpha a, \quad a_2=2^\beta b, \quad a_3=2^\gamma c,$$
where $a,b,c$ are positive odd integers.
The triviality of the intersection of 
$\Br(\ov X)[2]$ with the image of $\Br(X_k)$ in $\Br(\ov X)$
can be proved over a finite extension of $k$.
The field $K=k(\sqrt[4]{a}, \sqrt[4]{b}, \sqrt[4]{c})$ satisfies the assumptions of
Lemma \ref{good}, so we conclude the proof by replacing $k$ with $K$,
except in the case when $X_K$ is equivalent to $X_{1,2,8}\times_\Q K$ which is not covered by 
Lemma \ref{good}.

In the rest of the proof we deal with this exceptional case and assume $a_1=a$, $a_2=2 b$, $a_3=8 c$,
where $a,b,c$ are positive odd integers. 

\medskip

{\em The main case.} Let us first deal with the case $(a_1,a_2,a_3)=(1,2,8)$.
We observe that the product of the coefficients is a fourth power
in $\Q$, so the $\Ga_\Q$-modules $\Br(\ov X)$ and $\Br(\ov F)$ are isomorphic. 
We work over $K=k$ and apply Lemma \ref{new} (a). The action of 
$G_{a,K}=\Gal(k(\sqrt[4]{2})/k)$ on $\O/8$ factors through a group of order 2
whose generator acts as multiplication by $1+4i$.
By the standard formula for the first cohomology group of a cyclic group we obtain 
$\H^1(G_{a,K},\Delta_X)\simeq(\Z/2)^2$.
The injective map $\H^0(G_{a,K},\Br(\ov X)[2])=\O/2\to \H^1(G_{a,K},\Delta_X)$ 
is therefore an isomorphism.
By Lemma \ref{drei} the image of 
$\H^1(G_{a,K},\Pic(\ov X)^*)\to \H^1(G_{a,K},\Delta_X)$ has the same cardinality as the cokernel of 
$\H^0(G_{a,K},\Pic(\ov X)^*)\to \H^0(G_{a,K},\Delta_X)$. 
It can be checked on a computer (cf. Lemma~\ref{app-monster}) 
that this cardinality is 2.

To prove that $\left(\Br(X_k)/\Br_1(X_k)\right)\{2\}\cong\Z/2$ we show
that this group does not have an element of order 4.  
By Lemma \ref{drei} the image of the homomorphism
$\H^1(\widetilde G_{a,K},\Pic(\ov X)^*)\to \H^1(\widetilde G_{a,K},\Delta_X)$ 
has the same cardinality as the cokernel of 
$\H^0(\widetilde G_{a,K},\Pic(\ov X)^*)\to \H^0(\widetilde G_{a,K},\Delta_X)$. But this cokernel is the same as
the cokernel of $\H^0(G_{a,K},\Pic(\ov X)^*)\to \H^0(G_{a,K},\Delta_X)$, 
because the action of $\widetilde G_{a,K}$ factors through the action of $G_{a,K}$.
By Lemma \ref{new} (b) we conclude that the 2-primary torsion subgroup of $\Br(X_k)/\Br_1(X_k)$ is $\Z/2$.

\medskip

{\em The remaining cases.}
Having dealt with the main case we need to show that in all other cases
$\Br(\ov X)[2]$ intersects trivially with the image of $\Br(X_k)\to\Br(\ov X)$.
There is a prime $p$ dividing the positive odd integer $abc$. Let $L$ be the field
obtained by adjoining to $k$ the fourth roots of all the prime factors of $abc$ other than $p$.
Then $X_L$ is equivalent to a diagonal quartic surface of one of the following types.

$(p^n,2,8)$, where $n=1,2,3$. Let $K$ be the extension of $L$ obtained by adjoining 
the fourth root of $p/4$, $p/2$, $p/4$, respectively.
We note that $2$ is not a square in $K$ and $K_a=K(\sqrt[4]{2})$. 
As $X_K$ is equivalent to $X_{2,4,8}\times_\Q K$ considered in Lemma \ref{good},
we conclude by this lemma.

$(1,2p^m,8p^n)$, where $m>0$ or $n>0$. Let $K=L(\sqrt[4]{p/2})$. 
A tedious but easy verification shows that $X_K$ is equivalent to one of the
surfaces considered in Lemma \ref{good}, except in the case $(1,2p^2,8p^2)$.
In this case we let $M=L(\sqrt[4]{2p^2})$. Note that $2$ is a square in $M$,
so $8p^2$ is a fourth power in $M$, hence $M$ splits $X$. But $2$ is not a fourth power in $M$,
so we conclude by Theorem \ref{c1}.

$(p^2,2p^m,8p^n)$. Another tedious but easy verification shows that extending 
$L$ to $K=L(\sqrt[4]{p/2})$ reduces all cases, 
except those equivalent over $K$ to
$(p^2,2p^2,8)$, to the cases covered in Lemma \ref{good}. 

$(p,2p^m,8p^n)$. Let $K_0=L(\pi)$, where $\pi^2=p$, and let $K=K_0(\sqrt[4]{\pi/2})$.
We note that $2$ is not a square in $K$ and $K_a=K(\sqrt[4]{2})$. 
The surface $X_K$ is equivalent to one of the cases treated before, unless it is
equivalent to $(\pi^2,2\pi^2,8)$.

It remains to show that for the diagonal quartic surface $X$ with coefficients $(r^2,2r^2,8)$ over
a field $K\supset k$ such that $[K(\sqrt{2},\sqrt{r}):K]=4$ we have $\Br(X)\{2\}=\Br_1(X)\{2\}$.
The proof uses the following lemma.

\ble \label{C}
Let $A\times B$ be a product of two groups. Let $M$ be a $B$-module.
Consider $M$ as an $A\times B$-module with trivial action of $A$. Then
the inflation-restriction sequence with respect to the subgroup $A\subset A\times B$ is 
a split exact sequence, so we have a canonical isomorphism
$$\H^1(A\times B, M)\cong\H^1(A,M)^B\oplus\H^1(B,M),$$
where $\H^1(A,M)^B\cong\Hom(A,M^B)$.
\ele
{\em Proof.} We have $\H^1(A,M)\cong\Hom(A,M)$, and since $A$ and $B$ commute,
we have $\Hom(A,M)^B=\Hom(A,M^B)$.

The restriction map $\H^1(A\times B, M)\to\H^1(A,M)^B=\Hom(A,M^B)$
has a section given by the inflation map 
$\Hom(A,M^B)\to \H^1(A\times B,M)$. Thus $\H^1(A\times B, M)$ is the direct sum
of $\H^1(A,M)^B$ and $\H^1(B,M^A)=\H^1(B,M)$. $\Box$

\medskip

{\em End of proof of Theorem} \ref{main1}. 
The product of the coefficients $(r^2,2r^2,8)$ is a fourth power
in $K$, so the $\Gal(\ov K/K)$-modules $\Br(\ov X)$ and $\Br(\ov F)$ are isomorphic. 
We have $K_a=K(\sqrt[4]{2},\sqrt{r})$ and $G_{a,K}=G\times H$,
where $$G=\Gal(K(\sqrt[4]{2})/K)=\Gal(k(\sqrt[4]{2})/k)\simeq\Z/4, \quad\quad H=\Gal(K(\sqrt{r})/K)\simeq\Z/2,$$
and $H$ acts trivially on $\Br(\ov X)[16]$. The condition of Lemma \ref{new} is satisfied.

We can apply Lemma \ref{C} to $M=\Delta_X$, $A=H$, $B=G$.
Note that the inclusion $(\Delta_X)^G\to \Delta_X$ is the natural map $\O/4\to\O/8$;
it induces an isomorphism $\H^1(H,(\Delta_X)^G)\tilde\lra \H^1(H,\Delta_X)$ as both groups are $\O/2$.
By Lemma \ref{C} we obtain
$$\H^1(G\times H,\Delta_X)\cong\H^1(G,\Delta_X)\oplus \Hom(H,(\Delta_X)^G)
\cong\H^1(G,\Delta_X)\oplus \O/2.$$
The map $\partial_1\colon (\O/2)^G\to \H^1(G,\Delta_X)$ was already computed in the previous case
$(1,2,8)$, and was found to be an isomorphism. The inflation-restriction sequence with respect to the
subgroup $H\subset G\times H$ gives rise to the following commutative diagram, where the middle column
is a split exact sequence of Lemma \ref{C}:
$$\xymatrix{
\H^1(H,\Pic(\ov X)^*)\ar[r]&\O/2&\\
\H^1(G\times H,\Pic(\ov X)^*)\ar[r]\ar[u]^{\rm res}&\H^1(G,\Delta_X)\oplus 
\O/2\ar[u]^{\rm res}&(\O/2)^{G\times H}\ar[l]_{\quad\quad \partial_1}\\
&\H^1(G,\Delta_X)\ar@{^{(}->}[u]^{\rm inf}&(\O/2)^G\ar[l]_\cong\ar[u]^=
}$$
By Lemma \ref{new} (a), if the intersection of images of 
the maps in the middle row is zero, then $\Br(X_K)\{2\}=\Br_1(X_K)\{2\}$.
The commutativity of this diagram implies that ${\rm Im}(\partial_1)$ coincides
with the direct summand $\H^1(G,\Delta_X)$ of $\H^1(G,\Delta_X)\oplus \O/2$.
Thus we need to show that the image of $\H^1(G\times H,\Pic(\ov X)^*)$ intersects
$\H^1(G,\Delta_X)$ trivially. The cardinality
of the image of $\H^1(G\times H,\Pic(\ov X)^*)$ is $4$. This can be checked using 
the exact sequence in the proof of Lemma \ref{drei} using the fact that each of the groups
$\H^1(G\times H,\Pic(\ov X))$ and $\H^1(G\times H,\Pic(\ov X)^*)$ has $16$ elements.
Indeed, since $(\Delta_X)^{G\times H}\simeq (\Z/4)^2$, one needs to check that
the image of $(\Pic(\ov X)^*)^{G\times H}\to(\Delta_X)^{G\times H}$ is $4$.
Furthermore, the composition of the left upward
arrow and the top horizontal arrow is surjective; see Lemma~\ref{app-pseudomonster}
for all these facts. The commutativity of the diagram now implies
the desired result. This finishes the proof of Theorem \ref{main1}.
$\Box$

\subsection{The case of $\Q$}

The following is the main result of the paper.

\bthe \label{main2}
Let $X=X_a$ be a diagonal quartic surface with $a_1, a_2, a_3\in \Q^\times$.
Then $\left(\Br(X)/\Br_1(X)\right)\{2\}$ is $\Z/2$
if $X$ is equivalent to $X_{1,2,-2}$ or $X_{1,8,-8}$, and $0$ otherwise.
\ethe
{\em Proof.} By Theorem \ref{main1} we need to deal with the following cases:
\begin{equation}
(1,2,8), \quad (1,-2,-8), \quad (2,-2,-4), \quad (-2,4,8), \quad (1,2,-2), \quad (1,8,-8). \label{list}
\end{equation}
Each of these surfaces satisfies $\Q_a=\Q(i,\sqrt[4]{2})=k(\sqrt[4]{2})$, where $k=\Q(i)$.
By Lemma \ref{new}, in order to conclude that
$\Br(X)\{2\}=\Br_1(X)\{2\}$ it is enough to show that the images of
$\H^0(G_{a,\Q},\Br(\ov X)[2])$ and $\H^1(G_{a,\Q},\Pic(\ov X)^*)$ in $\H^1(G_{a,\Q},\Delta_X)$
intersect trivially. 

In the first three cases the product of coefficients is in $\Q^{\times4}$, so the $\Ga_\Q$-modules
$\Br(\ov X)$ and $\Br(\ov F)$ are isomorphic, hence
$\H^1(G_{a,\Q},\Delta_X)\simeq \H^1(G_{a,\Q},\Delta_F)$.
Proposition \ref{c2} (ii) implies that $\Br(\ov X)\{2\}^{\Ga_\Q}\simeq\Z/2\times\Z/4$,
hence $\H^0(G_{a,\Q},\Br(\ov X)[2])\cong\O/2$. 
A computer calculation (cf. Lemma~\ref{app-descent}) shows that $\H^1(G_{a,\Q},\Delta_X)\simeq(\Z/2)^4$ and in each case the image of $\H^1(G_{a,\Q},\Pic(\ov X)^*)$ has trivial intersection with the image of 
$\H^0(G_{a,\Q},\Br(\ov X)[2])$.

In the last three cases, where the product is in $-4\Q^{\times4}$, we have 
$\Br(\ov X)\{2\}^{\Ga_\Q}\simeq\Z/4$, hence $\H^0(G_{a,\Q},\Br(\ov X)[2])\cong\Z/2$. One checks that
in all three cases the group $\H^1(G_{a,\Q},\Delta_X)$ is $\Z/2$. 
Then a computer calculation (cf. Lemma~\ref{app-descent}) shows that in the case $(-2,4,8)$
the image of $\H^1(G_{a,\Q},\Pic(\ov X)^*)\to \H^1(G_{a,\Q},\Delta_X)$ is $0$. 
In contrast, in the last two cases
the map $\H^1(G_{a,\Q},\Pic(\ov X)^*)\to \H^1(G_{a,\Q},\Delta_X)$ turns out to be the surjective map
$\Z/4\to\Z/2$. Thus $\left(\Br(X)/\Br_1(X)\right)[2]\cong\Z/2$. 

It remains to prove that there is no element of order 4 in $\left(\Br(X)/\Br_1(X)\right)\{2\}$.
To do this, consider 
$\widetilde G_{a,\Q}=\Gal(k(\sqrt[8]{-2})/\Q)$. By Lemma \ref{new} it is enough to show
that the intersection of the image of $\H^0(\widetilde G_{a,\Q},\Br(\ov X)[4])$ with the image of
$\H^1(\widetilde G_{a,\Q},\Pic(\ov X)^*)$ in $\H^1(\widetilde G_{a,\Q},\Delta_X)$ is killed by 2. 
We have a commutative diagram
$$\xymatrix{
\H^1(\widetilde G_{a,\Q},\Pic(\ov X)^*)\ar[r]&\H^1(\widetilde G_{a,\Q},\Delta_X)\\
\H^1(G_{a,\Q},\Pic(\ov X)^*)\ar[r]\ar[u]_{\cong}&\H^1(G_{a,\Q},\Delta_X)\ar@{_{(}->}[u],
}$$
which shows that the image of the upper horizontal arrow is isomorphic to
the image of lower horizontal arrow, which is $\Z/2$. $\Box$

\bpr \label{comp}
Let $X$ be a diagonal quartic surface over $\Q$. Let $X'$ be
the diagonal quartic surface obtained by multiplying any two coefficients of $X$ by $-4$.
Then the $\Ga_\Q$-modules $\H_\et^2(\ov X,\Q_2)$ and $\H_\et^2(\ov X', \Q_2)$
are isomorphic, and the $\Ga_\Q$-modules $\H_\et^2(\ov X,\Z_\ell)$ and $\H_\et^2(\ov X', \Z_\ell)$ 
are isomorphic for all odd primes $\ell$. 
In particular, $|X(\F_q)|=|X'(\F_q)|$ for all prime powers $q$ coprime to $a_1a_2a_3$.
\epr
{\em Proof.} Twisting both surfaces we can assume that $X=F$ is the Fermat surface and 
$X'=F'=X_{(1,-4,-4)}$. Since the degree of the polarisation is $4$, it is enough to construct an
isomorphism between primitive cohomology $\Ga_\Q$-modules $P_\ell(F)$ and $P_\ell(F')$,
and also between $P(F)\otimes\Q_2$ and $P(F')\otimes\Q_2$.  Write $k=\Q(i)$.
Since $-4\in (k^\times)^4$, we have $F_k\cong F'_k$ and thus the representations 
of $\Ga_k$ are isomorphic. It remains to show that this isomorphism extends to complex conjugation.

Recall that the action of $\tau(1)$ on $P(F)$ fixes the generator: $\tau(1)(2\pi i e')=2\pi i e'$,
see (\ref{conj}). Thus for $P(F')$ we have $\tau(1)(2\pi i e')=u_2u_3(2\pi i e')$.
 In $\Z[x]/(x^3+x^2+x+1)$ we have $4(1-x)^{-1}=x^2 + 2x + 3$.
Hence $(1-x)(1-y)$ is an invertible element of $\Q_2[x,y]/(x^3+x^2+x+1,y^3+y^2+y+1)$
and of $\Z_\ell[x,y]/(x^3+x^2+x+1,y^3+y^2+y+1)$, where $\ell$ is odd.
It follows that multiplication by $(1-u_2)(1-u_3)$ defines isomorphisms
$P(F)\otimes\Q_2\cong P(F')\otimes\Q_2$ and $P_\ell(F)\cong P_\ell(F')$ where $\ell$ is odd.
Since $\Ga_k$ commutes with $(\mu_4)^3$, these isomorphisms respect the action of $\Ga_k$. As
one immediately checks, they also respect the action of $\tau(1)$. $\Box$

\appendix
\section{Computations for the main theorem}

This appendix collects details of computations 
used in the proofs of Theorems \ref{main1} and \ref{main2}. While these are not too
challenging to be unachievable by hand, it is more convenient to carry them out
with the help of a computer algebra system.
Magma \cite{Magma} code to this effect is available at 
\href{https://github.com/dgvirtz/diagonalquartics}{https://github.com/dgvirtz/diagonalquartics}
and from the auxiliary files accompanying this article 
{\tt arXiv:1905.11869} on the \href{https://arxiv.org}{arxiv.org} preprint server.

Consider the diagonal quartic surface $X_t$ with coefficients $t_1, t_2, t_3$
over the field $\Q(t_1,t_2,t_3)$, where $t_1, t_2, t_3$ are independent variables. Let
$$G_\eta=\Gal(\Q(i,\sqrt[4]{2},\sqrt[4]{t_1},\sqrt[4]{t_2},\sqrt[4]{t_3})/\Q(t_1,t_2,t_3)).$$ 
We denote the generators of $G_\eta$ by
$$\tau(i)=-i, \  \sigma(\sqrt[4]{2})=i\sqrt[4]{2}, \ 
  s_1(\sqrt[4]{t_1})= i\sqrt[4]{t_1}, \ 
  s_2(\sqrt[4]{t_2})= i\sqrt[4]{t_2}, \ 
  s_3(\sqrt[4]{t_3})= i\sqrt[4]{t_3},
$$
with the understanding that $\tau$ acts trivially on $\sqrt[4]{2},\sqrt[4]{t_1},\sqrt[4]{t_2},\sqrt[4]{t_3}$,
$\sigma$ acts trivially on $i, \sqrt[4]{t_1}, \sqrt[4]{t_2},  \sqrt[4]{t_3}$, and similarly for $s_1,s_2,s_3$.
A diagonal quartic surface $X_a$ over a number field $K$ gives rise 
to a subgroup $G_{a,K}\subset G_\eta$ which will be described explicitly in each case.

Fix fourth roots $\sqrt[4]{a_1},\sqrt[4]{a_2},\sqrt[4]{a_3}$ 
and a primitive eighth root of unity $\zeta$.
Recall that the abelian group $\Pic(\ov X)$ is generated by $48$ lines of the form
$$\Lambda_{h,i,j}^{m,n}:\quad
\zeta^m x_0=\sqrt[4]{a_h}x_h, \quad \zeta^n \sqrt[4]{a_i}x_i=\sqrt[4]{a_j}x_j,$$ 
where $(h,i,j)$ is a permutation of $(1,2,3)$ and $m,n\in\{1,3,5,7\}$. 
A basis $\B$ of $20$ lines, which Swinnerton-Dyer kindly provided to the authors, is given by
$$ \Lambda_{1,2,3}^{1,1},\Lambda_{1,2,3}^{1,3},\Lambda_{1,2,3}^{1,5},\Lambda_{1,2,3}^{1,7},\Lambda_{1,2,3}^{3,1},\Lambda_{1,2,3}^{3,3},\Lambda_{1,2,3}^{3,5},\Lambda_{1,2,3}^{5,1},\Lambda_{1,2,3}^{5,3},\Lambda_{1,2,3}^{5,5},$$
$$ \Lambda_{2,3,1}^{1,1},\Lambda_{2,3,1}^{1,3},\Lambda_{2,3,1}^{1,5},\Lambda_{2,3,1}^{3,1},\Lambda_{2,3,1}^{3,3},\Lambda_{2,3,1}^{3,5},$$
$$\Lambda_{3,1,2}^{1,1},\Lambda_{3,1,2}^{1,3},\Lambda_{3,1,2}^{1,5},\Lambda_{3,1,2}^{3,1}. $$
(A similar basis of $20$ lines appears in \cite{ShiShi}.)
The description of $\Pic(\ov X)$ in terms of lines can be matched with the description in the main text 
by first identifying the line $\Lambda_{1,2,3}^{1,1}$ 
with $\frac{1}{4}(L-c)$ 
and then identifying all other lines using the actions of $(\mu_4)^3$ and $\mathrm{S}_3$ on $\ov X$. 
The images of $w_1,w_2$ in $\Delta_X$ are expressed in the basis $\B$ as
$$\begin{array}{lllllllllllllllllllll}
\frac{1}{8}(1 & 1 & 0 & 0 & 1 & 4 & 7 & 0 & 7 & 3 & 6 & 6 & 4 & 2 & 2 & 4 & 0 & 4 & 0 & 4),\\
\frac{1}{8}(2 & 7 & 3 & 2 & 1 & 6 & 1 & 1 & 5 & 4 & 4 & 6 & 6 & 4 & 2 & 2 & 4 & 4 & 4 & 4).   
\end{array}$$
We write elements of $T(X_\C)$ and $\Delta_X\cong T(X_\C)/8$ in the basis $(w_1,w_2)$.
The 1-cocycles in group cohomology will be explicitly described as maps from $G_{a,K}$ to the coefficient module.

In order to assist an independent verification, we give intermediate steps and details of our computations.
Let $k=\Q(i)$ and let $G=\Gal(k(\sqrt[4]{2})/k)$.

\ble \label{app-surj}
Let $a\in \{(1,2,2),(1,2,4),(1,4,8)\}$.
 Let $K$ be a field extension of $k$ not containing $k(\sqrt{2})$.
 Then $G\cong G_{a,K}$ and the map $(\Pic(\ov X)^*)^G\to(\Delta_X)^G$ is surjective.
\ele
{\em Proof.}
In the first case, one has $G_{a,K}=\langle\sigma s_2 s_3\rangle$.
A computation of invariants gives $(\Pic(\ov X)^*)^G\simeq\Z^6$.
By Proposition \ref{c2}, one has $(\Delta_X)^G=\Delta_X[2]$; lifting
the generators of $\Delta_X[2]$ to $(\Pic(\ov X)^*)^G$ is easy. 

In the last two cases $G_{a,K}=\langle\sigma s_2 s_3^2\rangle$
and $G_{a,K}=\langle\sigma s_2^2 s_3^3\rangle$, respectively. 
A computation of invariants gives $(\Pic(\ov X)^*)^G\simeq\Z^3$. 
By Proposition \ref{c2}, the group $(\Delta_X)^G\simeq\Z/2$ is generated by $\frac{1}{2}(w_1+w_2)$ 
which lifts to $(\Pic(\ov X)^*)^G$. $\Box$

\ble\label{app-monster}
 Assume $K=k=\Q(i)$ and $a=(1,2,8)$. 
 Then $G\cong G_{a,k}$ and the cokernel of the map 
$(\Pic(\ov X)^*)^G\to(\Delta_X)^G$ has cardinality $2$.
\ele
{\em Proof.} One has
$G_{a,k}=\langle\sigma s_2 s_3^3\rangle$. Proposition \ref{c2} gives $(\Delta_X)^G=\Delta_X[4]$.
A computation of invariants shows that $(\Pic(\ov X)^*)^G\simeq\Z^6$ 
and that the image of this group in 
$\Delta_X[4]$ is generated by the elements $\frac{1}{4}(w_1\pm w_2)$. $\Box$

\ble\label{app-pseudomonster}
Let $K$ be an extension of $k=\Q(i)$ and let $a=(r^2,2r^2,8)$, where $r\in K^\times$ is such that
$[K(\sqrt{2},\sqrt{r}):K]=4$. We have $G_{a,K}\cong G\times H$, where
$H=\Gal(K(\sqrt{r})/K)$. 
 The image of $(\Pic(\ov X)^*)^{G\times H}\to(\Delta_X)^{G\times H}$ has cardinality $4$, 
$$|\H^1(G\times H, \Pic(\ov X)^*)|=|\H^1(G\times H, \Pic(\ov X))|=16$$ and the composition
of the following maps is surjective:
 $$\H^1(G\times H, \Pic(\ov X)^*)\stackrel{\rm res}\lra\H^1(H, \Pic(\ov X)^*)\lra 
\H^1(H, \Delta_X)=\Delta_X[2].$$
\ele
{\em Proof.}
 One has $G_{a,K}=\langle\sigma s_2 s_3^2, s_1^2s_2^2\rangle$.
 A computation gives $(\Pic(\ov X)^*)^{G\times H}\simeq \Z^4$. 
The four generators of $(\Pic(\ov X)^*)^{G\times H}$ lie in $\Pic(\ov X)$ only 
after doubling and, furthermore, they are not all in the same equivalence class modulo $\Pic(\ov X)$. 
Thus their images span $\Delta_X[2]$.

 The first cohomology of $\Pic(\ov X)$ was already computed in \cite[Table B]{Bright}, case B37. 
Another computation shows that $\H^1(G\times H, \Pic(\ov X)^*)$ 
is generated by three 1-cocycles $c_1$, $c_2$ and $c_3$ of orders $2$, $2$ and $4$ 
such that $c_2$ and $c_3$ map to generators of $\H^1(H, \Delta_X)\cong\Delta_X[2]$. $\Box$

\medskip

Lemmas A.1--A.3 concern $\Pic(\ov X), \Pic(\ov X)^*$ and $\Delta_X$. 
In particular, no computational input about the transcendental lattice and its interplay with the 
algebraic cycles is necessary. 
Thus the verification of the computational facts up to this point can be achieved using 
solely the description of $\Pic(\ov X)$ in terms of lines. This is not so for the next lemma used in the proof of Theorem \ref{main2}.

\ble\label{app-descent}
 Assume that $K=\Q$ and $a$ is in the list $(\ref{list})$. Then the intersection of images
of the following maps
$$\H^1(G_{a,\Q}, \Pic(\ov X)^*)\stackrel{\rho}\lra\H^1(G_{a,\Q}, \Delta_X)\stackrel{\partial_1}
\longleftarrow \Br(\ov X)[2]^{G_{a,\Q}}$$ 
is trivial in the first four cases and is $\Z/2$ in the last two cases of $(\ref{list})$.
\ele
{\em Proof.}
In all cases $G_{a,\Q}$ is generated by two elements $g_1,g_2$, which we specify. 
We describe $1$-cocycles by the images of $g_1$ and $g_2$.

\smallskip

\noindent\textbf{Case $a=(1,2,8)$.} We have $g_1=\sigma s_2s_3^3$ and $g_2=\tau$.
A computation gives $\H^1(G_{a,\Q}, \Pic(\ov X)^*)\simeq(\Z/2)^2$.
The group $\H^1(G_{a,\Q}, \Delta_X)\simeq(\Z/2)^4$ is generated by
$$\begin{array}{lcrr}
g_1&\mapsto& (0, & 2) \\
g_2&\mapsto& 0
  \end{array} \quad\quad
\begin{array}{lcrr}
g_1&\mapsto& (-4, & 0) \\
g_2&\mapsto& 0
  \end{array} \quad\quad
\begin{array}{lcrr}
g_1&\mapsto& 0 \\
g_2&\mapsto& (-4, & 0)
  \end{array} \quad\quad
\begin{array}{lcrr}
g_1&\mapsto& (-2, & 0) \\
g_2&\mapsto& (0, & 1).
  \end{array}
$$
In this basis of $\H^1(G_{a,\Q}, \Delta_X)$, we have
$$\Im(\partial_1)=\langle \begin{matrix} (1, & 1, & 0, & 0)\end{matrix},\begin{matrix} (0, & 0, & 0, & 1)\end{matrix}\rangle
\text{ and }
\Im(\rho)=\langle\begin{matrix} (0, & 1, & 0, & 0)\end{matrix}\rangle.$$

\noindent\textbf{Case $a=(1,-2,-8)$.} We have $g_1=\sigma s_2s_3^3$ and $g_2=s_2s_3\tau$.
A computation gives $\H^1(G_{a,\Q}, \Pic(\ov X)^*)\simeq(\Z/2)^2$.
The group $\H^1(G_{a,\Q}, \Delta_X)\simeq(\Z/2)^4$ is generated by
$$\begin{array}{lcrr}
g_1&\mapsto& (4, & 4) \\
g_2&\mapsto& 0
  \end{array}\quad\quad
\begin{array}{lcrr}
g_1&\mapsto& (-2, & 0) \\
g_2&\mapsto& 0
  \end{array}
\quad\quad\begin{array}{lcrr}
g_1&\mapsto& 0 \\
g_2&\mapsto& (4, & 4)
  \end{array}
\quad\quad\begin{array}{lcrr}
g_1&\mapsto& (2, & 2) \\
g_2&\mapsto& (-1, & 0).
  \end{array}
$$
In this basis, $\Im(\partial_1)$ is as in the previous case 
and $\Im(\rho)=\langle \begin{matrix} (1, & 0, & 1, & 0) \end{matrix}\rangle$.

\noindent\textbf{Case $a=(2,-2,-4)$.} We have $g_1=\sigma s_1s_2^3$ and $g_2=s_2s_3\tau$.
A computation gives an isomorphism $\H^1(G_{a,\Q}, \Pic(\ov X)^*)\simeq(\Z/2)^3$.
The group $\H^1(G_{a,\Q}, \Delta_X)\simeq(\Z/2)^4$ is generated by the same cocycles as for $a=(1,-2,-8)$.
In this basis, $\Im(\partial_1)$ is as in the previous case and 
$\Im(\rho)=\langle \begin{matrix} (1, & 0, & 0, & 0) \end{matrix},\begin{matrix} (0, & 0, & 1, & 0) \end{matrix}\rangle.$

\noindent\textbf{Case $a=(-2,4,8)$.} We have $g_1=\sigma s_1^3s_2^2s_3^3$ and $g_2=s_1\tau$.
Then $\H^1(G_a, \Pic(\ov X)^*)\simeq(\Z/2)^2$ is generated by two cocycles
which map to $0$ in $\H^1(G_a, \Delta_X)$.

\smallskip

\noindent\textbf{Case} $a=(1,2,-2)$ (respectively, $a=(1,8,-8)$). 
We have $g_1=\sigma s_2s_3^3$ (respectively, $g_1=\sigma s_2^3s_3$) and $g_2=s_3\tau$.
A computation shows that $\H^1(G_a, \Pic(\ov X)^*)\simeq\Z/4$
maps surjectively onto $\H^1(G_a, \Delta_X)\simeq\Z/2$, which is generated by
the cocycle sending $g_1$ to $(-2,2)$ and $g_2$ to $0$. $\Box$

\section{Computations for the supplement}

Let $a=(1,2,-2)$ or $(1,8,-8)$. 
In view of Theorems \ref{main1} and \ref{main2} we can apply
Corollary \ref{1.5} (b) with $n=2$. Thus we need
to compute the first hypercohomology group of 
the finite group $G_{a,K}$ with coefficients in the complex $\Pic(\ov X)^*\to T(X_\C)/16$,
where $K=\Q$ and $\Q(i)$. 

Hypercohomology functionality is currently not available in Magma and had to be added.
Our implementation constructs a part of the Cartan--Eilenberg resolution of the above complex and computes its first total cohomology. More precisely, let $\ZZ^1(M)$ be the group of $1$-cocycles of a 
$G_{a,K}$-module $M$. Consider the double complex
$$\xymatrix{
0&&\\
\ZZ^1(\Pic(\ov X)^*)\ar[r]^{h_{01}} \ar[u]&\ZZ^1(T(X_\C)/16) &\\
\Pic(\ov X)^*\ar[r]^{h_{00}} \ar[u]^{v_{00}}&T(X_\C)/16\ar[r]\ar[u]^{v_{10}} & 0
}$$
where $h_{01}$ is the map induced by $h_{00}$ and the vertical maps are the differentials of the standard resolutions sending $m$ to $gm-m$. Details of the calculations are as follows. (Matrices representing linear maps are understood to act on the right.)

\noindent\textbf{Case $K=\Q(i)$.}
One computes that $\ZZ^1(\Pic(\ov X)^*)$ is freely generated by $14$ cocycles
and $\ZZ^1(T(X_\C)/16)$ is generated by two cocycles which map 
the generator $\sigma s_2s_3^3$ of $G_{a,K}$ to $(4,4)$ and to $(-4,0)$.
Using $(w_1,w_2)$ and the above basis for $\ZZ^1(T(X_\C)/16)$, one checks that $v_{10}$ is given by
$$\begin{bmatrix}
-1 & -3\\
-2 & 9
  \end{bmatrix}.$$
Similarly, one easily finds matrices representing $h_{01}$ and $v_{00}$.
A computation gives that $\Ker(h_{01}-v_{10})\cong(\Z/4)^2\oplus\Z^{14}$ and, after choosing generators $x_1,x_2$ of the torsion subgroup and linearly independent $y_1,\dots,y_{14}$ suitably in this kernel, we 
find that $\Im(v_{00}\oplus h_{00})$ is generated by
$y_1+y_2+2y_3,2y_1+x_1,2y_2+x_1,x_2,y_4,\dots,y_{14}$.
Thus the first hypercohomology is $(\Z/4)^2$.

\noindent\textbf{Case $K=\Q$.} Recall that $G_{a,\Q}$ is generated by $g_1=\sigma s_2s_3^3$ 
(or $\sigma s_2^3s_3$) and $g_2=s_3\tau$.
One computes that $\ZZ^1(\Pic(\ov X)^*)$ is freely generated by $17$ cocycles
and $\ZZ^1(T(X_\C)/16)$ is generated by the two cocycles
$$\begin{array}{lcrr}
g_1 &\mapsto &(4 & 0)\\
g_2 &\mapsto &(-1 & 1)
  \end{array},
\quad \quad\begin{array}{lcrr}
g_1 &\mapsto & (-4 & 4)\\
g_2 &\mapsto & 0
  \end{array}.
$$
Using $(w_1,w_2)$ and the above basis for $\ZZ^1(T(X_\C)/16)$, one checks that $v_{10}$ is represented by
$$\begin{bmatrix}
1 & -1\\
-1 & 2
  \end{bmatrix}.$$
Similarly, one easily finds matrices representing $h_{01}$ and $v_{00}$.
A computation gives that $\Ker(h_{01}-v_{10})\cong\Z/4\oplus\Z^{17}$ and, after choosing $x$ of order $4$ and linearly independent $y_1,\dots,y_{17}$ suitably, $\Im(v_{00}\oplus\phi)$ is generated by 
$2x,x+4y_1,y_2,\dots,y_{17}$.
Thus the first hypercohomology is $\Z/8$.

\bigskip

\noindent Department of Mathematics, University College London, WC1E 6BT England, U.K.

\medskip

\noindent d.gvirtz@ucl.ac.uk

\bigskip

\noindent Department of Mathematics, South Kensington Campus,
Imperial College London, SW7 2BZ England, U.K. -- and --
Institute for the Information Transmission Problems,
Russian Academy of Sciences, 19 Bolshoi Karetnyi, Moscow, 127994
Russia

\medskip

\noindent a.skorobogatov@imperial.ac.uk

\end{document}